\documentclass[preprint,12pt]{elsarticle}

\usepackage{amsmath,amssymb,amsthm}
\usepackage{algorithm,algorithmic}
\usepackage[labelformat=simple]{subfig}

\usepackage{epsfig}
\usepackage{mathptmx,bm}
\usepackage{multirow}

\newtheorem{theorem}{Theorem}[section]

\newtheorem{remark}[theorem]{Remark}

\newcommand{\set}[2]{\left\{{#1}\,:~{#2}\right\}}
\newcommand {\average}[1] {\mbox{$\left\{\!\!\left\{ #1 \right\}\!\!\right\}$}}
\newcommand {\jump}[1] {\mbox{$\left[\!\left[ #1 \right]\!\right]$}}
\newcommand {\la}[1] {\left\langle {#1} \right\rangle }

\journal{Computers and Mathematics with Applications}

\begin{document}

\begin{frontmatter}

\title{Reduced Order Optimal Control of the Convective FitzHugh-Nagumo Equations}

\author[iam]{B\"ulent Karas\"{o}zen \corref{cor1}}
\ead{bulent@metu.edu.tr}

\author[sin]{Murat Uzunca}
\ead{muzunca@sinop.edu.tr}

\author[bal]{Tu\u{g}ba K\"u\c{c}\"ukseyhan }
\ead{kucukseyhan@balikesir.edu.tr}

\cortext[cor1]{Corresponding author}

\address[iam]{Institute of Applied Mathematics \& Department of Mathematics, Middle East Technical University, 06800 Ankara, Turkey}
\address[sin]{Department of Mathematics, Sinop University, 57000 Sinop, Turkey}
\address[bal]{Department of Mathematics, Bal{\i}kesir  University,  10145 Bal{\i}kesir,  Turkey}

\begin{abstract}
In this paper, we compare three model order reduction methods: the proper orthogonal decomposition (POD), discrete empirical interpolation method (DEIM) and dynamic mode decomposition (DMD) for the optimal control of the convective  FitzHugh-Nagumo (FHN) equations. The convective  FHN equations consists of the semi-linear activator and the linear inhibitor equations, modeling blood coagulation in moving excitable media. The semilinear activator equation leads to a non-convex optimal control problem (OCP).  The most commonly used method in reduced optimal control is POD. We use DEIM and DMD  to approximate efficiently the nonlinear terms in reduced order models. We compare the accuracy and computational times of three reduced-order optimal control solutions with the full order discontinuous Galerkin finite element solution of the convection dominated FHN equations with terminal controls.  Numerical results show that POD is the most accurate whereas POD--DMD is the fastest.
\end{abstract}


\end{frontmatter}

\section{Introduction}
\label{introduction}

Optimal control problems (OCPs) governed by semilinear partial differential equations (PDEs) with wave-type solutions  have been recently investigated for the Schl\"ogl equation with traveling waves fronts  and for the Nagumo equation with spiral waves \cite{buchholz13ocs}, for the  FitzHugh-Nagumo (FHN) equations with spiral and traveling waves \cite{Casas13,Casas15soasa,Ryll16}, lambda-omega systems with spiral waves \cite{borzi06doc}, Allee equation \cite{Trelat18}, optimal harvesting \cite{Yilmaz17} and  the convective FHN equations \cite{Uzunca17} with traveling wave solutions. PDE-constrained optimization is also applied to parameter identification in pattern formation. Optimal parameters leading to Turing patterns in semi-linear reaction-diffusion equations like the Schnakenberg and Gierer-Meinhardt equations are identified in \cite{Garvie10,Garvie14,Stoll16}. Controlling traveling wave solutions or patterns of semi-linear PDEs is computationally  challenging. For stable solutions of these OCPs, at each iteration step of the nonlinear optimization algorithms, large scale linear systems should be solved. Efficient and stable solutions of this kind of PDE constrained optimization problems are the focus of the above mentioned studies.

For space discretization of PDEs, commonly used methods are finite differences, finite elements, finite volume, and spectral methods.  One of the most stable and accurate methods is the discontinuous Galerkin (dG) method for convection dominated PDEs like the FHN equations.  In this paper, we use the symmetric interior penalty Galerkin (SIPG) method for space discretization \cite{Arnold02,riviere08dgm}, which is the most relevant and popular dG method. The dG methods possess better properties for convective problems compared with continuous finite element methods, but they require a much larger number of degrees of freedom. Recently OCPs governed by linear stationary and time-dependent convection-diffusion-reaction equations \cite{Akman14,Akman14a,Yucel14} and semilinear steady state equations \cite{Yucel15dgm} are solved efficiently by the several dG methods. OCPs with PDE constraints are discretized usually following two approaches: the {\em discretize-then-optimize\/} approach, where the optimality conditions are imposed in the discrete setting, or  the {\em optimize-then-discretize} approach in which the optimality conditions are derived in variational form and then discretization holds.  There is no preferred approach  \cite{Benner14tip,Gubisch17}, here we follow the {\em optimize-then-discretize} approach: the optimality conditions are discretized in space by the SIPG method and the backward Euler method in time. We employ the projected nonlinear conjugate gradient (CG) method  \cite{Hager06acgd} for solving the resulting nonlinear discrete optimization problem.

The optimal control of instationary PDEs requires numerous evaluations of the optimality system, which lead to large scale optimization problems. In the last decade, different types of reduced order modeling (ROM) methods were developed to approximate these by low order or surrogate models. The resulting small optimization problems can be solved efficiently with less computational effort. For an overview of ROMs for OCPs for linear and semilinear PDEs, we refer the reader to \cite{Benner14tip,Gubisch17}. Despite its heuristic nature, currently, proper orthogonal decomposition (POD) is the preferred model reduction technique for linear and semilinear  OCPs. The POD basis functions are constructed from the snapshots, solutions of space-time discretized PDE at pre-specified time-instances, then the reduced order OCP  is solved applying the Galerkin projection.   ROMs using the POD-Galerkin projection for OCPs for linear PDE constraints are investigated in \cite{Hinze08,Kunisch08pod,Troltzsch08} and for semi-linear PDE constraints in   \cite{Kammann13,Studinger13,Kunisch15}. With increasing number of POD basis functions, more accurate solutions are obtained. But the success of POD depends on the type of the problem; for semi-linear PDEs with traveling or spiral wave type solutions like the Schl\"ogl and FHN equations \cite{buchholz13ocs,Ryll16,Ryll14pod} require more POD modes than for the OCPs with semi-linear parabolic equations \cite{Kammann13}. Because the reduced nonlinear term still depends on the dimension of the full order model (FOM), i.e., high dimensional finite element discretized model, hyper reduction techniques are developed.  The empirical interpolation method (EIM) \cite{Barrault04} and the discrete empirical interpolation method (DEIM) \cite{chaturantabut10nmr} are commonly used as hyperreduction techniques to reduce the computational cost for evaluation of the nonlinear terms.  More recently, as an alternative to the DEIM, the dynamic mode decomposition (DMD) \cite{Alla16} is used to approximate the nonlinear terms of the PDEs. DMD was first introduced by Schmid \cite{Schmid10dmd}, Rowley \cite{Rowley12} and it is based on the linear approximation of the infinite dimensional Koopman operator \cite{Koopman31,Kutz16}.

We compare the accuracy and computational efficiency of  POD, POD with DEIM and POD with DMD for the optimal control of the convective FHN equations modeling blood coagulation and bio-reactors \cite{Ermakova09opo,Ermakova05}. The FHN equations are the simplest and most widely used model for describing the complex spatio-temporal behavior of traveling waves in excited media.  In contrast to the classical FHN  equations \cite{Fitzhugh61,Nagumo62} with a semi-linear PDE and with a linear ordinary differential equation (ODE), the convective FHN equations consists of a  semi-linear PDE with the monotone cubic nonlinear term for the activator and a linear PDE for the inhibitor, modeling excitable systems in moving media. For semi-linear OCPs like the Burger's equation, quadratic nonlinear terms  can be reduced to the linear terms \cite{Baumann18}. But for higher order polynomial or general nonlinearities, hyperreduction techniques should be used.  As far as we know, reduced order OCPs using POD--DEIM and POD--DMD are not yet investigated for semilinear PDE constraint OCPs.  Among these three methods, POD--DMD is the fastest, because after collecting snapshots, the non-linearity disappears and the reduced model becomes a coupled system of linear ODEs. Also, the reduced non-convex optimization problem for POD and POD--DEIM becomes convex for POD--DMD.  The POD is the most accurate but, the slowest. The POD--DMD is less accurate than the POD but fastest.

In the sequel, we describe the OCP constrained by the convective FHN equations in Section~\ref{sec:ocp}. The fully discrete optimal control system is given in Section~\ref{ocp_disc}, based on the first order optimality conditions. In Section~\ref{rom}, the reduced order OCPs (POD, POD--DEIM, and POD--DMD) are constructed. In Section~\ref{numeric}, we compare the three ROMs with respect to accuracy and computational efficiency for a test problem with terminal control. The paper ends with concluding remarks in Section~\ref{conc}.

\section{Optimal control problem}
\label{sec:ocp}

We consider the following OCP:
\begin{equation} \label{ocp}
\begin{aligned}
\min \limits_{f \in \mathcal{F}_{ad}} J(u,v,f)
     = & \; \frac{1}{2}  \int \limits_{\Omega} \big( u(x,T)-u_T(x) \big)^2 \; dx
       + \frac{1}{2}  \int \limits_{\Omega} \big( v(x,T)-v_T(x) \big)^2 \; dx  \\
     & \; +  \frac{\upsilon}{2}   \int \limits_{0}^{T} \int \limits_{\Omega}  f(x,t)^2 \; dx\;dt,
\end{aligned}
\end{equation}
governed by the convective FHN equations
\begin{equation}\label{state}
\begin{aligned}
u_t(x,t) - d_u \Delta u(x,t) + \bm{b}(x) \cdot \nabla u(x,t) + g(u(x,t)) + v(x,t)   &= f(x,t) & \text{ in } Q,  \\
v_t(x,t) - d_v \Delta v(x,t) + \bm{b}(x) \cdot \nabla v(x,t) + \epsilon (v(x,t) - c_3 u(x,t)) &= 0 & \text{ in } Q,\\
\partial_{\emph{n}} u(x,t) = 0, \quad  \partial_{\emph{n}} v(x,t) =  0  && \text{on   } \Sigma^{N},\\
u(x,t) =  u_D(x,t), \quad  v(x,t) = v_D(x,t) && \text{on  } \Sigma^{D}.
\end{aligned}
\end{equation}
where $d_u$, $d_v$ are the diffusion coefficients, $c_3$, $\epsilon$ are real constants, $x=(x_1,x_2)^T$ is the spatial element, $u_T(x)$ and $v_T(x)$ are given desired terminal states, and $u_D(x,t)$ and $v_D(x,t)$ are given functions.
The control function $f(x,t)$  satisfy the point-wise box constraints
\begin{equation*}
f \in \mathcal{F}_{ad} := \{ f \in L^{\infty}(Q): \; f_l \leq f(x,t) \leq f_r \;\; \hbox{for   a.e   } (x,t) \in Q \},
\end{equation*}
for some real numbers $-\infty < f_l < f_r < +\infty $.
The cubic polynomial $g(u)$ is given by
\begin{equation*}
g(u)=c_{1} u(u-c_{2})(u-1),
\end{equation*}
where $c_{i}$,  $i=1,2$, are non-negative real numbers. The non-linearity $g(u)$ is monostable for $0 < c_1< 20$ and $c_2=0.02$  \cite{Ermakova09opo}, whereas for the Schl\"ogl equation \cite{buchholz13ocs},  the classical FHN equations \cite{Casas13}, and the diffusive FHN equations \cite{Karasozen15} the bistable cubic non-linearity is  considered.  The velocity field $\bm{b}(x)= (b_1(x), b_2(x))^T$ is divergence-free, and for the maximum wave speed $b_{\max}$, it is given by the parabolic profile along the $x_1$-direction
\begin{equation*}
b_1(x) = ax_2(H-x_2), \quad b_{max}=\frac{1}{4}aH^2, \quad a >0 , \quad b_2(x) = 0.
\end{equation*}

In \eqref{ocp} and \eqref{state}, $T>0$ and $Q:= \Omega \times (0,T]$ denote the final time and the time-space domain, respectively.  The space domain $\Omega=(0,L) \times (0,H)$ is bounded, Lipschitz  in $\mathbb{R}^2$, and $\Sigma=\partial\Omega \times (0,T]$ is the lateral surface. The  Dirichlet and Neumann boundaries are $\Sigma^D := \Gamma_D \times (0,T]$ and $\Sigma^N := \Gamma_N \times (0,T]$, where $\Gamma_D=\{x=\{0,L\},\; 0 < y < H \}$ and  $\Gamma_N =\{0 < x < L,\; y =\{0,H\} \}$. The outward unit normal vector is  denoted by $\bm{n}(x)$ and the outward normal derivative on $\partial \Omega$ is defined as
 $\partial_{\emph{n}}u=\nabla u\cdot \bm{n}$.

The convective FHN equations \eqref{state} consists of an activator equation for $u(x,t)$ and an inhibitor equation for $v(x,t)$ in an excitable medium. In blood coagulation process, the activator variable $u(x,t)$  describes the concentrations of thrombin in the excitation, and $v(x,t)$ the inhibition
of this excitation and recovery of the medium variable, which is activated factor XI \cite{Ermakova05}. The complex process of coagulation consists of cascadic  enzymatic reactions and feedback loops \cite{Lonabov05}. These reactions enable generation of auto-catalytic thrombin outside of the
damaged region.  The most important property of the blood coagulation process is the formation of auto-waves with the velocity independent
initial conditions  \cite{Ermakova09opo,Ermakova05,Lonabov05}.
The flow propagates  within the impermeable channel walls (Neumann boundary conditions).  Different types of standing and triggering waves occur depending on the constants of the convective FHN equations \cite{Ermakova09opo}. The waves become more curved when the flow velocity $b_{max}$ increases. For a detailed discussion of the complex spatio-temporal wave phenomena occurring in the uncontrolled convective FHN equations, we refer the reader to \cite{Ermakova09opo,Ataullakhanov}.

The waves are initiated at $t=0$ inside an initial
narrow excitation region $\{ x_a \le x_1 \le x_b, \; 0 \le x_2 \le H\}$
\begin{equation} \label{init}
u_0(x) = \left\{
             \begin{array}{ll}
               \bar{u}, & \hbox{if} \;\; x_a \leq x_1 \leq x_b,  \\
               0, & \hbox{otherwise},
             \end{array}
           \right.
 \quad
v_0(x) = 0.
\end{equation}
The interval from $x_a$ to $x_b$ is the damaged region
of the vessel wall. In  blood coagulation, for fixed thrombin concentration
the initial condition \eqref{init}   mimics the activation process \cite{Ermakova09opo,Ermakova05,Lonabov05}. In this study, we take a rectangle as $x_a=0$, $x_b=0.1$, and we set $\bar{u}=0.1$.

Optimal control aims that the  activator $u(x,t)$ and the inhibitor $v(x,t)$ follow desired terminal states $u_T(x)$ and $v_T(x)$ as close as possible in the $L_2$ norm by minimizing the cost functional $J(u,v,f)$. The inhibitor $v(x,t)$ has only some auxiliary character as for the FHN equations in \cite{Casas13,Casas15soasa,Ryll16}, therefore only the activator $u(x,t)$  is to be controlled by $f(x,t)$, which seems to be experimentally feasible \cite{Vilas08}. The parameter $\upsilon> 0$ in the optimal cost functional $J(u,v,f)$ denotes the penalization or the Tikhonov regularization parameter.

\begin{remark}
We have considered the OCP \eqref{ocp} with terminal observation functionals to show the applicability of reduced order modeling. Reduced solutions of OCPs with observation functionals, distributed in space–time for the FHN equations \cite{Uzunca17} can also be handled using the same approach.
\end{remark}

\section{Optimality system and discretization}
\label{ocp_disc}

The local optimal solutions are constructed using the first-order necessary optimality conditions. In \cite{Uzunca17}, we consider a more general OCP constrained by the FHN equations with distributed, terminal and sparse controls. Here we consider the OCP \eqref{ocp} only with the terminal controls and without sparse controls.
As in \cite{Uzunca17}, we follow the {\em optimize-then-discretize} approach, first by deriving the first order optimality conditions in the variational form and then discretize the OCP. The  second order necessary and sufficient optimality conditions for semi-linear parabolic equations are investigated in \cite{Casas13, Troltzsch10}.
The first-order necessary optimality conditions \cite[Chapter 1]{Troltzsch10},  lead to the following adjoint system
\begin{equation}\label{adjoint}
\begin{aligned}
-p_t(x,t) - d_u\Delta p(x,t) - \bm{b}(x)\cdot\nabla p(x,t)  + g'(u)p(x,t) - \epsilon c_3q(x,t) &= 0,  \\
-q_t(x,t) - d_v\Delta q(x,t) - \bm{b}(x)\cdot\nabla q(x,t) + \epsilon q(x,t) + p(x,t) &= 0,\\
\end{aligned}
\end{equation}
with the mixed boundary and terminal time conditions
\begin{equation*}
\begin{aligned}
d_u \partial_{\emph{n}} p(x,t) + (\bm{b} \cdot \bm{n}) p(x,t) = 0, \quad & d_v \partial_{\emph{n}} q(x,t) + (\bm{b} \cdot \bm{n})q(x,t)  =  0, & \text{on  } \Sigma^N, \\
p(x,t) =  0, \quad & q(x,t)  =  0, &  \text{on  } \Sigma^D,  \\
p(x,T) =   u(x, T)-u_T(x), \quad & q(x,T)  = v(x, T)-v_T(x), &  \text{in  } \Omega.
\end{aligned}
\end{equation*}
The flow in the adjoint system \eqref{adjoint} is in the opposite direction of the flow in the state system \eqref{state}.
In view of the box constraints, the optimal solution $f(x,t)$ has to satisfy the following variational inequality
\begin{align*}
\int \limits_{0}^{T} \int \limits_{\Omega} \big( p(x,,t) + \upsilon f(x,t)  \big) \big(\nu(x,t) -f(x,t)\big) \;dx \;dt \geq 0, \qquad \forall \nu \in \mathcal{F}_{ad},
\end{align*}
which leads to the point-wise projection formula
\begin{equation*}
f(x,t) = \mathbb{P}_{[f_l,f_r]}\left( -\frac{1}{\upsilon}p(x,t) \right),
\end{equation*}
with the projection operator $\mathbb{P}_{[a,b]}(\phi)= \max \{a,\min \{\phi,b\}\}$.

\subsection{Space discretization of the optimality system}
\label{disc_state}

The optimality system consisting of the state equation \eqref{state} and the adjoint equation \eqref{adjoint} is discretized in space using the SIPG method \cite{riviere08dgm}. We denote the family of meshes consisting of non-overlapping triangular elements $K$  by $\{\mathcal{T}_h\}_h$. The diameter of an element $K$ and the length  of an edge $E$ are denoted by $h_K$ and $h_E$, respectively. We assume
that the mesh is regular in the following sense: for different triangles $K_i,\; K_j  \in \{\mathcal{T}_h\}_h, i \neq j$, the intersection $K_i \cap K_j$ is either empty or a vertex or an edge, i.e., hanging nodes are not allowed.

The discrete test, state, and control spaces are defined  by the  space of piecewise  discontinuous  finite element functions
\begin{equation*}
W_{h} = \set{w \in L^2(\Omega)}{ w\mid_{K}\in \mathbb{P}^l(K) \; , \; \forall K \in \mathcal{T}_h}.
\end{equation*}
Here, the set of  polynomials on $K \in \mathcal{T}_h$ of degree at most $l$ is denoted by $\mathbb{P}^l(K)$. In numerical examples, we use linear discontinuous finite elements, i.e., $l=1$.

The set $\mathcal{E}_h$ of all edges are split into the sets of interior, Dirichlet and Neumann boundary edges, denoted respectively by $\mathcal{E}_h^0,\, \mathcal{E}_h^D,\, \mathcal{E}_h^N$, so that $\mathcal{E}_h = \mathcal{E}_h^0 \cup \mathcal{E}_h^D \cup \mathcal{E}_h^N$. The inflow and outflow boundaries are defined by
$$
\Gamma^- = \set{x \in \partial \Omega}{ \bm{b}(x) \cdot \bm{n}(x) < 0}, \quad
\Gamma^+ = \partial \Omega \setminus  \Gamma^-.
$$
Similarly, the inflow and outflow boundaries of an element $K \in \mathcal{T}_h$ are defined by
\begin{equation*}
\partial K^- =\set{x \in \partial K}{\bm{b}(x) \cdot \bm{n}_{K}(x) <0}, \quad \partial K^{+} = \partial K \setminus \partial K^{-},
\end{equation*}
where $\bm{n}_{K}(x)$ is the exterior unit normal vector on the boundary $\partial K$ of an element $K$.
Let the common edge for two elements $K$ and $K^e$ be $E$. Then there are two traces of a function $w\in W_h$ along the edge $E$,  denoted by $w|_E$ from the interior of  $K$ and $w^e|_E$ from the interior of  $K^e$. Accordingly, the jump and average of $w$ across the edge $E$ are defined by
\begin{equation*}
\jump{w} =w|_E\bm{n}_{K} + w^e|_E\bm{n}_{K^e}, \quad
\average{w}=\frac{1}{2}\big( w|_E + w^e|_E \big).
\end{equation*}
Similarly,  the jump and average of a vector-valued function are given by
\begin{equation*}
\jump{\nabla w} =\nabla w|_E \cdot \bm{n}_{K}+\nabla w^e|_E \cdot \bm{n}_{K^e}, \quad
\average{\nabla w}=\frac{1}{2}\big(\nabla w|_E+\nabla w^e|_E \big).
\end{equation*}
We set $\jump{w}=w|_E\bm{n}$ and $ \average{\nabla w}= w|_E$ on a boundary edge $E \in K \cap \partial\Omega$.  Then, the state equations \eqref{state} and the adjoint equations \eqref{adjoint} are discretized using the SIPG method leading to the system: $\forall w \in W_h$ and a.e. $t\in (0,T]$
\begin{equation}\label{sipgstate}
\begin{aligned}
\la{\partial_t u_h, w } + a_{h,u}(u_h,w) + g_h(u_h,w) + \la{v_h,w} &= \ell_{h,u}(t,w) + \la{f_h,w}, \\
\la{u_h(x ,0), w(x)} &= \la{u_0(x ),w(x)}, \\
\la{ \partial_t v_h, w } + a_{h,v}(v_h,w) + \epsilon \la{v_h,w} - \epsilon c_3\la{u_h,w} &= \ell_{h,v}(t,w), \\
\la{v_h(x,0), w(x)}&= \la{v_0(x),w(x)},
\end{aligned}
\end{equation}

\begin{equation}\label{sipgadj}
\begin{aligned}
-\la{ \partial_t p_h, w } + a_{h,p}(p_h,w) + \la{g'(u_h)p_h,w}-\epsilon c_3\la{q_h,w} &= 0, \\
\la{p_h(x,T),w(x)} &=  \la{ u_h(x, T)-u_T(x) ,w(x) },  \\
-\la{ \partial_t q_h, w } + a_{h,q}(q_h,w) + \epsilon \la{q_h,w}  + \la{p_h,w} &= 0, \\
 \la{q_h(x,T),w(x)} &=  \la{ v_h(x, T)-v_T(x) ,w(x) },
\end{aligned}
\end{equation}
where $\la{\cdot , \cdot}$ denotes the usual $L^2$-inner product. The bilinear forms $a_{h,u}$ and $a_{h,v}$ in the discrete state system \eqref{sipgstate} are given for $i=u,v$ by
\begin{equation*}
\begin{aligned}
a_{h,i}(\tilde{w},w)=& \sum \limits_{K \in \mathcal{T}_h} \int \limits_{K} d_i\nabla \tilde{w}  \cdot  \nabla w \; dx - \sum \limits_{E \in \mathcal{E}^0_h \cup \mathcal{E}^D_h} \int \limits_E \Big( \average{ d_i\nabla  \tilde{w}} \cdot \jump{w} + \average{ d_i\nabla w} \cdot \jump{\tilde{w}} \Big)  \; ds   \\
& + \sum \limits_{E \in \mathcal{E}^0_h \cup \mathcal{E}^D_h} \frac{d_i\gamma}{h_E} \int \limits_E \jump{\tilde{w}} \cdot \jump{w} \; ds + \sum \limits_{K \in \mathcal{T}_h} \int \limits_{K} \bm{b} \cdot \nabla \tilde{w} w  \; dx \\
&  + \sum \limits_{K \in \mathcal{T}_h}\; \int \limits_{\partial K^{-} \backslash \partial \Omega} \bm{b} \cdot \bm{n} (\tilde{w}^e-\tilde{w})w \; ds - \sum \limits_{K \in \mathcal{T}_h} \; \int \limits_{\partial K^{-} \cap \Gamma^{-}} \bm{b} \cdot \bm{n} \tilde{w} w  \; ds,
\end{aligned}
\end{equation*}
for any $\tilde{w},\;w \in W_h$, and $\gamma$ is called the penalty parameter \cite{riviere08dgm}.
For $i=u,v$, the linear right-hand side vectors, and the nonlinear form are defined by
\begin{align*}
\ell_{h,i}(t,w) &= \sum \limits_{E \in \mathcal{E}^D_h} \int \limits_E  i_{D}(x,t) \Big(\frac{d_i\gamma}{h_E} \bm{n} \cdot \jump{w} - \average{d_i \nabla w}  \Big)\; ds
 - \sum \limits_{K \in \mathcal{T}_h} \; \int \limits_{\partial K^{-} \cap \Gamma^{-}} \bm{b} \cdot \bm{n} \; i_{D}(x,t) w  \; ds,\\
g_h(u,w) &= \sum \limits_{K \in \mathcal{T}_h} \int \limits_{K} g(u)  w  \; dx.
\end{align*}
The bilinear forms $a_{h,p}$ and $a_{h,q}$ in the adjoint system \eqref{sipgadj} are defined similar to the bilinear forms $a_{h,u}$ and $a_{h,v}$ for states, respectively, but with the negative velocity field. They contain, respectively, the additional terms
$$
\sum \limits_{E \in \mathcal{E}^N_h}  \int \limits_E \big( \bm{b} \cdot \bm{n} \big) p_h w\; ds \quad \hbox{and} \quad \sum \limits_{E \in \mathcal{E}^N_h}  \int \limits_E \big( \bm{b} \cdot \bm{n} \big) q_h w\; ds,
$$
which come from the mixed boundary conditions.
The SIPG discrete states and the control in \eqref{sipgstate} are given in the form
\begin{equation*}
u_h(x,t) = \sum \limits_{i=1}^{N_K} \sum \limits_{j=1}^{N_{l}} u_{j}^{\,i}(t) \phi_{j}^{\,i}(x), \;
v_h(x,t) = \sum \limits_{i=1}^{N_K} \sum \limits_{j=1}^{N_{l}} v_{j}^{\,i}(t) \phi_{j}^{\,i}(x), \;
f_h(x,t) = \sum \limits_{i=1}^{N_K} \sum \limits_{j=1}^{N_{l}} f_{j}^{\,i}(t) \phi_{j}^{\,i}(x).
\end{equation*}
In a similar way, the discrete desired terminal states (just $L^2$-projections) are given by
\begin{equation*}
u_{h,T}(x) = \sum \limits_{i=1}^{N_K} \sum \limits_{j=1}^{N_{l}} u_{T,j}^{\,i} \phi_{j}^{\,i}(x), \quad
v_{h,T}(x) = \sum \limits_{i=1}^{N_K} \sum \limits_{j=1}^{N_{l}} v_{T,j}^{\,i} \phi_{j}^{\,i}(x).
\end{equation*}
The number of (triangular) elements are denoted by $N_K$, and the local dimension of each element is given by $N_{l}$. $\phi_{j}^{\,i}$ denotes the $j$-th finite element basis function defined on the $i$-th triangle. Setting dG degrees of freedom $N:=N_K\times N_l$, the corresponding time dependent unknown coefficient vectors can be written as
\begin{align*}
\bm{u}(t) &= (u_{1}^{\,1}(t), \ldots, u^{\,1}_{N_l}(t), \ldots, u^{\,N_K}_{1}(t), \ldots, u^{\,N_K}_{N_l}(t))\in\mathbb{R}^N, \\
\bm{v}(t) &= (v_{1}^{\,1}(t), \ldots, v^{\,1}_{N_l}(t), \ldots, v^{\,N_K}_{1}(t), \ldots, v^{\,N_K}_{N_l}(t))\in\mathbb{R}^N, \\
\bm{f}(t) &= (f_{1}^{\,1}(t), \ldots, f^{\,1}_{N_l}(t), \ldots, f^{\,N_K}_{1}(t), \ldots, f^{\,N_K}_{N_l}(t))\in\mathbb{R}^N.
\end{align*}
On the other hand, the known coefficient vectors of the discrete desired terminal states are given by
\begin{align*}
\bm{u}_T &= (u_{T,1}^{\,1}, \ldots, u^{T,1}_{N_l}, \ldots, u^{T,N_K}_{1}, \ldots, u^{T,N_K}_{N_l})\in\mathbb{R}^N, \\
\bm{v}_T &= (v_{T,1}^{\,1}, \ldots, v^{T,1}_{N_l}, \ldots, v^{T,N_K}_{1}, \ldots, v^{T,N_K}_{N_l})\in\mathbb{R}^N,
\end{align*}
which are computed from the projection identities
\begin{equation*}
\la{u_{h,T}(x),\phi_i(x)} = \la{u_T(x),\phi_i(x)} \; , \quad  \la{v_{h,T}(x),\phi_i(x)} = \la{v_T(x),\phi_i(x)}\; , \quad i=1,\ldots N.
\end{equation*}

Then, the SIPG semi-discretized system of the state equations \eqref{sipgstate} leads to the following system of ODEs
\begin{equation}\label{full_order_ode_state}
\begin{aligned}
 \bm{M} \bm{u}_t +  {\bm S}_u\bm{u} +  \bm{g}(\bm{u}) + \bm{M}\bm{v}  &=  { \bm \ell}_u + \bm{M} \bm{f}, \\
 \bm{M} \bm{v}_t +  {\bm S}_v\bm{v} + \epsilon \bm{M} \bm{v} -\epsilon c_3 \bm{M} \bm{u} &=   \bm{\ell}_v,
\end{aligned}
\end{equation}
and the semi-discrete objective functional takes form
\begin{align*}
J_h(\bm{u},\bm{v},\bm{f}) = & \; \frac{1}{2}(\bm{u}-\bm{u}_T)^T\bm{M}(\bm{u}-\bm{u}_T) + \frac{1}{2}(\bm{v}-\bm{v}_T)^T\bm{M}(\bm{v}-\bm{v}_T) \\
& \; + \frac{\upsilon}{2}\int_0^T\bm{f}^T\bm{M}\bm{f}dt.
\end{align*}
Here $\bm{M}$ denotes the mass matrix, and  $\bm{S}_u$ and $\bm{S}_v$ are the stiffness matrices corresponding to the bilinear forms $a_{h,u}$ and $a_{h,v}$, respectively. The nonlinear vector $\bm{g}(\bm{u})$ corresponds to  the nonlinear form $g_h$, and the time-dependent vectors ${ \bm \ell}_u:={ \bm \ell}_u(t)$ and ${ \bm \ell}_v:={ \bm \ell}_v(t)$ are the vectors corresponding to the linear forms $\ell_{h,u}$ and $\ell_{h,v}$, respectively. By a similar setting, the SIPG semi-discretized system of the adjoint equations \eqref{sipgadj} are given as the following system of ODEs
\begin{equation}\label{full_order_ode_adj}
\begin{aligned}
 -\bm{M} \bm{p}_t +  {\bm S}_p\bm{p} + \bm{R}_g(\bm{u})\bm{p} -\epsilon c_3 \bm{M}\bm{q}  &=  0, \\
 -\bm{M} \bm{q}_t +  {\bm S}_q\bm{q} + \epsilon \bm{M} \bm{q} + \bm{M} \bm{p} &=  0,
\end{aligned}
\end{equation}
where the state $u(x,t)$ dependent matrix $\bm{R}_g(\bm{u})$ is related to the term $\la{g'(u_h)p_h,w}$ in the adjoint system \eqref{sipgadj}, and $\bm{S}_p$ and $\bm{S}_q$ are the stiffness matrices corresponding to the bilinear forms $a_{h,p}$ and $a_{h,q}$, respectively.

\subsection{Fully discrete optimality system}
\label{fullydisc_state}

The time interval $[0,T]$ is partitioned uniformly as $0=t_0 < t_1 < \cdots < t_{N_T}=T$, with the step-size $\Delta t = T/N_T$. For $n=0,1,\ldots , N_T$, we denote by $\bm{u}_n$, $\bm{v}_n$, $\bm{f}_n$, $\bm{p}_n$ and $\bm{q}_n$ the approximate solution vectors to the semi-discrete solutions $\bm{u}(t)$, $\bm{v}(t)$, $\bm{f}(t)$, $\bm{p}(t)$ and $\bm{q}(t)$ at the time instance $t=t_n$, respectively. Also, we set $\bm{\ell}_u^n:=\bm{\ell}_u(t_n)$ and $\bm{\ell}_v^n:=\bm{\ell}_v(t_n)$. Then, applying the  backward Euler method to the semi-discrete state equation \eqref{full_order_ode_state} and to the semi-discrete adjoint equation \eqref{full_order_ode_adj},  we obtain the following fully discrete state system
\begin{equation}\label{fos}
\begin{aligned}
 \frac{1} {\Delta t} \bm{M} \big( \bm{u}_n - \bm{u}_{n-1} \big)  +   {\bm S}_u \bm{u}_{n} +  \bm{g}(\bm{u}_{n}) + \bm{M} \bm{v}_{n}  =  \bm{\ell}_u^n + \bm{M} \bm{f}_n, \\
 \frac{1} {\Delta t } \bm{M} \big( \bm{v}_n - \bm{v}_{n-1} \big) +  {\bm S}_v \bm{v}_n   + \epsilon \bm{M} \bm{v}_n  - \epsilon c_3 \bm{M} \bm{u}_{n} =  \bm{\ell}_v^n,\\
(\bm{M} \bm{u}_0)_i =  \la{u_0(x),\phi_i}, \quad (\bm{M} \bm{v}_0)_i = \la{v_0(x),\phi_i},
\end{aligned}
\end{equation}
for $n=1,2,\ldots, N_T$, and the following fully discrete adjoint system
\begin{equation}\label{foa}
\begin{aligned}
 \frac{1} {\Delta t } \bm{M} \big( \bm{p}_{n-1} - \bm{p}_{n} \big)  + \bm{S}_p \bm{p}_{n-1}  + \bm{R}_g(\bm{u}_{n}) \bm{p}_{n-1} - \epsilon c_3 \bm{M}\bm{q}_{n-1}  &=  0, \\
 \frac{1} {\Delta t} \bm{M} \big( \bm{q}_{n-1} - \bm{q}_{n} \big) +  \bm{S}_q \bm{q}_{n-1}  + \epsilon \bm{M} \bm{q}_{n-1} + \bm{M} \bm{p}_{n-1} &=  0,\\
\bm{p}_{N_T} = \bm{u}_{N_T}  -  \bm{u}_T, \quad
\bm{q}_{N_T} = \bm{v}_{N_T}  -  \bm{v}_T,&
\end{aligned}
\end{equation}
for $n= N_T, \ldots,2,1$, where $\bm{u}_T$ and $\bm{v}_T$ are the known coefficient vectors of the desired terminal states. We note that the fully discrete state system \eqref{fos} is solved forward in time, whereas the fully discrete adjoint system \eqref{foa} is solved backward in time. The equations \eqref{fos} and \eqref{foa} are semi-explicit in time. Because the mass matrix ${\mathbf M}$ is symmetric positive definite and therefore invertible, the discrete solutions of \eqref{fos} and \eqref{foa} exist and are unique.
The resulting semilinear discrete OCP can be solved by several optimization algorithms.
We use the  projected  nonlinear CG method  \cite{Hager06acgd}, which is applied for solving the OCP governed by Schl\"ogl and FHN equations in \cite{buchholz13ocs,Casas13, Casas15soasa,Ryll16} and by the convective FHN equations in \cite{Uzunca17}. The details of the implementation of the projected nonlinear CG algorithm can be found in these papers.

\section{Reduced order optimal control}
\label{rom}

We construct ROMs for the state equation \eqref{state} and the adjoint equation \eqref{adjoint}. For the construction of the POD modes to the state variables $u$ and $v$, the POD method uses the snapshot matrices generated by the coefficient vectors $\bm{u}$ and $\bm{v}$  of the optimal state solution vectors from the full order optimal control problem, and similarly, the snapshot matrix generated by the coefficient vector $\bm{f}$ of the optimal control solution is used to construct the POD modes to the control variable $f$. For the construction of the POD modes to the adjoint variables $p$ and $q$, we use the same POD modes generated by the snapshots of the states $u$ and $v$, motivated by the error analysis in  \cite{Gubisch17}. This approach might be not the best option, but  the construction of different POD modes using the snapshot matrices of the adjoint states does not improve the suboptimal solutions much \cite{Kammann13} and  requires more computational work.
It is not known a-priorly how far the optimal solution
of the reduced POD problem is from the exact one unless the snapshots are generating a sufficiently large state space. Different approaches are discussed  in \cite{Benner14tip} for the proper choice of the POD basis in the reduced OCP, which affects the accuracy of the reduced solutions. Here we construct the POD basis for the reduced OCP from the snapshots of controls of the FOM.

\subsection{POD Galerkin discretization}
\label{podGalerkin}

The reduced order system for the state equations \eqref{state} of lower dimensions $k_u$ and $k_v$ is formed by the Galerkin projection of the equations onto  the subspaces
\begin{equation*}
W_{u,h}^{r}=\text{span} \{ \psi_{u,1}, \ldots, \psi_{u,k_u} \}\subset W_h, \quad W_{v,h}^{r}=\text{span} \{ \psi_{v,1}, \ldots, \psi_{v,k_v}\}\subset W_h.
\end{equation*}
For any time $t$, the reduced order states $u_{h}^r(x,t)\in W_{u,h}^{r}$ and $v_{h}^r(x,t)\in W_{v,h}^{r}$ stand for the approximations of the full order states $u_{h}(x,t)$ and $v_{h}(x,t)$, respectively, and they are given by
\begin{equation*}
u_h(x,t)\approx u_{h}^r(x,t) = \sum_{i=1}^{k_u} u_i^r(t) \psi_{u,i}(x), \quad v_h(x,t)\approx v_h^r(x,t)=\sum_{i=1}^{k_v} v_i^r(t) \psi_{v,i}(x),
\end{equation*}
where $\bm{u}^r(t):= ( u^r_1(t),\ldots, u^r_{k_u}(t) )^T$ and $\bm{v}^r(t):= ( v^r_1(t),\ldots, v^r_{k_v}(t) )^T$ are the reduced coefficient vectors, and $\{\psi_{u,i} \}_{i=1}^{k_u}$ and $\{\psi_{v,i} \}_{i=1}^{k_v}$ are the $L^2$-orthogonal reduced basis functions. Similarly, the full order control $f_{h}(x,t)$ is approximated by the reduced  order control $f_{h}^r(x,t)$ from the $k_f$ dimensional subspace $W_{f,h}^{r}$
\begin{equation*}
f_h(x,t)\approx f_{h}^r(x,t) = \sum_{i=1}^{k_f} f_i^r(t) \psi_{f,i}(x), \quad W_{f,h}^{r}=\text{span} \{ \psi_{f,1}, \ldots, \psi_{f,k_f} \}\subset W_h,
\end{equation*}
where $\bm{f}^r(t):= ( f^r_1(t),\ldots, f^r_{k_f}(t) )^T$ is the reduced coefficient vector for the control, and $\{\psi_{f,i} \}_{i=1}^{k_f}$ are the $L^2$-orthogonal reduced basis functions. All the reduced basis functions are linear combinations of dG basis functions
\begin{equation*}
\psi_{u,i} = \sum_{j=1}^{N} \Psi_{u,j,i} \phi_j(x), \quad \psi_{v,i} = \sum_{j=1}^{N} \Psi_{v,j,i} \phi_j(x), \quad \psi_{f,i} = \sum_{j=1}^{N} \Psi_{f,j,i} \phi_j(x),
\end{equation*}
where $\Psi_{\cdot,\cdot,i}$ are the coefficient vectors of the $i$-th reduced basis functions. For $w\in\{u,v,f\}$,
the following POD matrices are constructed
\begin{equation*}
\bm{\Psi}_w := [ \Psi_{w,1} \cdots \Psi_{w,k_w}] \in\mathbb{R}^{N\times k_w}, \quad \Psi_{w,i} = \left( \Psi_{w,1,i}, \ldots ,\Psi_{w,N,i} \right)^T \in\mathbb{R}^N,
\end{equation*}
where the columns $\Psi_{w,i}$ are the POD modes.
The POD modes are computed through the snapshot matrices $\mathcal{U}=[\bm{u}_1,\ldots ,\bm{u}_{N_T}]$, $\mathcal{V}=[\bm{v}_1,\ldots ,\bm{v}_{N_T}]$ and $\mathcal{F}=[\bm{f}_1,\ldots ,\bm{f}_{N_T}]$ in $\mathbb{R}^{N\times N_T}$ \cite{Karasozen15,Kunisch01}, respectively, where the columns of the snapshot matrices are simply the coefficient vectors of the full order states and control from the FOM \eqref{fos} at the time instances $t_i$,  $i=1,\ldots , N_T$. Then, for $w\in\{u,v,f\}$, the reduced basis functions $\{\psi_{w,i}\}$, $i=1,2,\ldots , k_w$, are given by the solution of the following minimization problem
\begin{align*}
\min_{\psi_{w,1},\ldots ,\psi_{w,k_w}} \frac{1}{N_T}\sum_{j=1}^{N_T} \left\| \bm{w}_j - \sum_{i=1}^{k_w} \la{\bm{w}_j,\psi_{w,i}}\psi_{w,i}\right\|_{L^2(\Omega)}^2 \\
\text{subject to } \la{\psi_{w,i},\psi_{w,j}} = \Psi_{w,\cdot ,i}^T{\bm M}\Psi_{w,\cdot ,j}=\delta_{ij} \; , \; 1\leq i,j\leq k_w,
\end{align*}
where $\delta_{ij}$ is the Kronecker delta. Again for $w\in\{u,v,f\}$, and for $\mathcal{W}\in\{\mathcal{U},\mathcal{V},\mathcal{F}\}$, the above minimization problem is equivalent to the eigenvalue problems
\begin{equation}\label{eg1}
\mathcal{W}\mathcal{W}^T{\bm M}\Psi_{w,\cdot ,i}=\sigma_{w,i}^2\Psi_{w,\cdot ,i} \; , \quad i=1,2,\ldots ,k_w,
\end{equation}
for the coefficient vectors $\Psi_{w,\cdot ,i}$ of the reduced basis functions $\psi_{w,i}$. Setting $\widehat{\mathcal{W}}=R\mathcal{W}$ ($R^T$ is the Cholesky factor of the mass matrix ${\bm M}$), we obtain the equivalent formulation of  \eqref{eg1} as
\begin{equation*}
\widehat{\mathcal{W}}\widehat{\mathcal{W}}^T\widehat{\Psi}_{w,\cdot ,i}=\sigma_{w,i}^2\widehat{\Psi}_{w,\cdot ,i} \; , \quad \quad i=1,2,\ldots ,k_w,
\end{equation*}
where $\widehat{\Psi}_{w,\cdot ,i}=R\Psi_{w,\cdot ,i}$.
Because the singular value decomposition (SVD) is more stable and efficient than the eigenvalue decomposition, we apply the SVD to compute the first $k_w$ left singular vectors, $\widehat{\Psi}_{w,\cdot ,i}:=\zeta_{w,i}$, of the  matrix $\widehat{\mathcal{W}}$
\begin{equation*}
\widehat{\mathcal{W}} = \zeta_{w}\Sigma_w \beta_w^T \; ,
\end{equation*}
where the diagonal matrix $\Sigma_w$ contains the singular values $\sigma_{w,i}$, and $\zeta_{w,i}$ are the columns of the orthogonal matrix $\zeta_{w}$. Using  the fact that $\widehat{\Psi}_{w,\cdot ,i}=R\Psi_{w,\cdot ,i}$, the coefficient vectors $\Psi_{w,\cdot ,i}$ of the reduced basis functions are computed as
$$
\Psi_{w,\cdot ,i}=R^{-1}\widehat{\Psi}_{w,\cdot ,i} \; , \quad i=1,2,\ldots k_w.
$$
By the above settings, the following relations can be found between the coefficient vectors $\bm{u}$, $\bm{v}$ and $\bm{f}$ of the FOM solutions, and the reduced coefficient vectors $\bm{u}^r$, $\bm{v}^r$ and $\bm{f}^r$
\begin{eqnarray*}
\bm{u} \approx \bm{\Psi}_u \bm{u}^r, & \bm{v} \approx \bm{\Psi}_v \bm{v}^r, & \bm{f} \approx \bm{\Psi}_f \bm{f}^r,\\
\bm{u}^r \approx \bm{\Psi}_u^T{\bm M} \bm{u}, & \bm{v}^r \approx \bm{\Psi}_v^T{\bm M} \bm{v}, & \bm{f}^r \approx \bm{\Psi}_f^T{\bm M} \bm{f}.
\end{eqnarray*}
Then the reduced system for the state equations takes the form
\begin{equation}\label{staterom}
\begin{aligned}
  \frac{d}{dt} \bm{u}^r +  \bm{S}^r_u \bm{u}^r   + \bm{\Psi}_u^T \bm{g}(\bm{\Psi}_u \bm{u}^r)
  + \bm{M}_{u,v}^r \bm{v}^r  &=  \bm{\Psi}_u^T \bm{\ell}_u + \bm{M}_{u,f}^r \bm{f}^r, \\
 \frac{d}{dt} \bm{v}^r +  \bm{S}^r_v \bm{v}^r  + \epsilon \bm{v}^r
- \epsilon c_3 \bm{M}^r_{v,u} \bm{u}^r &=  \bm{\Psi}_v^T \bm{\ell}_v,
\end{aligned}
\end{equation}
and using the $L^2$-orthogonality of the reduced basis functions, the semi-discrete reduced order objective functional takes form
\begin{align*}
J_h^r(\bm{u}^r,\bm{v}^r,\bm{f}^r) = & \; \frac{1}{2}(\bm{u}^r-\bm{u}^r_T)^T(\bm{u}^r-\bm{u}^r_T) + \frac{1}{2}(\bm{v}^r-\bm{v}^r_T)^T(\bm{v}^r-\bm{v}^r_T) \\
& \; + \frac{\upsilon}{2}\int_0^T(\bm{f}^r)^T\bm{f}^rdt,
\end{align*}
which is absolutely cheaper to calculate, comparing to the full order objective functional $J_h$.
Using the same POD modes generated by the snapshot matrix of the states, we can obtain the reduced adjoint system in a similar way
\begin{equation}\label{romadjoint}
\begin{aligned}
  -\frac{d}{dt} \bm{p}^r +  \bm{S}_p^{\,r} \bm{p}^r + \bm{R}_g^{\,r}(\bm{u}^{r}) \bm{p}^r
  -\epsilon c_3 \bm{M}_{u,v} \bm{q}^r  &=  0, \\
 -\frac{d}{dt} \bm{q}^r +  \bm{S}_q^{\,r} \bm{q}^r
 + \epsilon \bm{q}^r  +  \bm{M}_{v,u}^{\,r} \bm{p}^r  &= 0.
\end{aligned}
\end{equation}
In the ROMs \eqref{staterom} and \eqref{romadjoint}, the reduced matrices are defined by
\begin{align*}
& \bm{S}^r_u=\bm{\Psi}_u^T \bm{S}_u\bm{\Psi}_u, \quad \bm{S}^r_v=\bm{\Psi}_v^T \bm{S}_v\bm{\Psi}_v, \quad \bm{S}^r_p=\bm{\Psi}_u^T \bm{S}_p\bm{\Psi}_u, \quad \bm{S}^r_q=\bm{\Psi}_v^T \bm{S}_q\bm{\Psi}_v  \\
& \bm{M}_{u,v}^{\,r}=\bm{\Psi}_u^T \bm{M} \bm{\Psi}_v, \quad \bm{M}_{u,f}^{\,r}=\bm{\Psi}_u^T \bm{M} \bm{\Psi}_f, \quad \bm{M}_{v,u}^{\,r}=\bm{\Psi}_v^T \bm{M} \bm{\Psi}_u, \\
& \bm{R}_g^{\,r}(\bm{u}^{r})=\bm{\Psi}_u^T \bm{R}_g(\bm{\Psi}_u\bm{u}^{r})\bm{\Psi}_u.
\end{align*}
The sub-optimality systems  \eqref{staterom} and \eqref{romadjoint} are solved by  backward Euler method in time.

\subsection{Discrete empirical interpolation method}
\label{deim}

In this part, we rewrite the reduced state equation \eqref{staterom} resulting from the application of the DEIM \cite{chaturantabut10nmr} to the cubic nonlinear term $ \bm{g}(\bm{\Psi}_u \bm{u}^r)$.
The DEIM aims to approximate the nonlinear vector $ \bm{g}(\bm{\Psi}_u \bm{u}^r)$ by projecting it onto a subspace of the space generated by the non-linear functions and spanned by a basis of dimension $m\ll N$
\begin{equation}\label{podG}
 \bm{g}(\bm{\Psi}_u \bm{u}^r) \approx \bm{W}\bm{s}(t),
\end{equation}
where $\bm{s}(t)\in\mathbb{R}^m$ is the corresponding coefficient vector, and the DEIM basis matrix $\bm{W}\in\mathbb{R}^{N\times m}$ is calculated by the application of the POD to the nonlinear snapshot matrix $\mathcal{G}=[\bm{g}(\bm{u}_1) \cdots  \bm{g}(\bm{u}_{N_T})]\in\mathbb{R}^{N\times N_T}$. Since the system \eqref{podG} is overdetermined, we find a projection matrix $\bm{P}$ such that $\bm{P}=[\bm{e}_{\mathfrak{p}_1} \cdots  \bm{e}_{\mathfrak{p}_m}]\in\mathbb{R}^{N\times m}$ where $\bm{e}_{\mathfrak{p}_i}$ is the $i$-th column of the identity matrix $\bm{I}\in\mathbb{R}^{N\times N}$. After elimination of the coefficient vector $\bm{s}(t)$ in \eqref{podG}, we arrive at the following ROM for the state equation with the DEIM approximation as follows
\begin{equation}\label{statedeim}
\begin{aligned}
  \frac{d}{dt} \bm{u}^r +  \bm{S}^r_u \bm{u}^r   + \bm{Q}\bm{g}^r(\bm{u}^r)
  + \bm{M}_{u,v}^r \bm{v}^r  &=  \bm{\Psi}_u^T \bm{\ell}_u + \bm{M}_{u,f}^r \bm{f}^r, \\
 \frac{d}{dt} \bm{v}^r +  \bm{S}^r_v \bm{v}^r  + \epsilon \bm{v}^r
- \epsilon c_3 \bm{M}^r_{v,u} \bm{u}^r &=  \bm{\Psi}_v^T \bm{\ell}_v,
\end{aligned}
\end{equation}
where the matrix $\bm{Q}:=\bm{\Psi}_u^T \bm{W}(\bm{P}^T\bm{W})^{-1}$ is computed once in the off-line stage. The computation of the  reduced nonlinear vector $\bm{g}^r(\bm{u}_r):= \bm{P}^T \bm{g}(\bm{\Psi}_u \bm{u}^r)\in\mathbb{R}^m$ and its Jacobian matrix require only $m\ll N$ and $m\times N_l$ integral evaluations, respectively, while they need $N$ and $N\times N_l$ integral evaluations when DEIM approximation is not used.

\subsection{Dynamic mode decomposition}
\label{dmd}

The DMD extracts dynamically relevant spatio-temporal information content from  numerical or experimental data sets \cite{Kutz16}.  Without explicit knowledge of the dynamical system, the DMD algorithm determines eigenvalues, eigenmodes, and spatial structures for each mode.  The snapshot data
is decomposed in spatio-temporal modes both by POD and DMD that correlates the spatial features of the data. Additionally, DMD
associates the snapshot data to the  temporal Fourier modes.
DMD is a special case of the Koopman operator  \cite{Koopman31} approximating nonlinear systems via an associated infinite dimensional system. The connection between the DMD and Koopman operator was  established  in \cite{Rowley12,Schmid10dmd} . The Koopman operator ${\mathcal K}$ acts on a set of scalar observable functions $g:{\mathcal M} \rightarrow {\mathbb C}$
\begin{equation*}
{\mathcal K} g({\bm y}) = g({\bm N}({\bm y})) ,
\end{equation*}
for the nonlinear dynamical system
\begin{equation*}
\frac{d{\bm y}}{dt} = {\bm N}({\bm y}),
\end{equation*}
where ${\bm y}\in {\mathcal M}$, and ${\mathcal M}$ is an n-dimensional manifold. The DMD determines the Koopman eigenvalues and  modes directly from the data, when linear observables are considered as state space, $g({\bm y}) = {\bm y}$, directly from the data.

We consider a snapshot matrix and a time shifted version of the snapshot matrix
\begin{equation}\label{dmdsnap}
{\bm Y}= [ {\bm y}(t_0) \cdots  {\bm y}(t_{m-1}) ], \quad
{\bm Y}'= [ {\bm y}(t_1)  \cdots  {\bm y}(t_m) ].
\end{equation}
The DMD involves the decomposition of the unknown best-fit linear operator ${\mathbf A}\in \mathbb{R}^{n \times n}$ relating the matrices above
\begin{equation*}
{\mathbf Y}' = {\mathbf A} {\mathbf Y}.
\end{equation*}
The DMD modes are computed by the exact DMD algorithm \cite{Tu14}.

\begin{algorithm}[t]
 \caption{Exact DMD Algorithm\label{exactDMD_alg}}
 \begin{algorithmic}
 \STATE ${\mathcal G }$ and ${\mathcal G }'$ are the snapshot matices
 \STATE Computation of  SVD of ${\mathcal G }$,  ${\mathcal G } = U\Sigma V^{*}$.
 \STATE Define $\tilde{\bm A}_{\mathcal G } = U^{*} {\mathcal G }' V \Sigma^{-1}$.
 \STATE Computation of the eigenvalues and eigenvectors   $(\lambda_j,w_j)$ from $\tilde{\bm A}_{\mathcal G } w_j=  \lambda_jw_j$, $j=1,\ldots ,\widetilde{m}$.
 \STATE Computation of DMD modes $\bm{\Psi}^{\text{DMD}}_j:={\mathcal G }'V\Sigma^{-1}w_j$, $j=1,\ldots, \widetilde{m}$.
 \end{algorithmic}
 \end{algorithm}

In \cite{Alla16}, the DMD is proposed for producing low-rank
approximations of the nonlinearities in PDEs as an alternative to DEIM. Here, after building the POD basis functions of rank $k_u$,
we collect snapshots  ${\mathcal G}, {\mathcal G}'$  for the  non-linear term ${\bm g}(\bm{u})$ from the FOM \eqref{full_order_ode_state} at $N_{T+1}$ time instances and divide them into sets as in \eqref{dmdsnap}
$$
{\mathcal G} = [ {\bm g}(\bm{u}_0) \cdots {\bm g}(\bm{u}_{N_{T-1}}) ]\in\mathbb{R}^{N\times N_T}, \quad   {\mathcal G}'   = [ {\bm g}(\bm{u}_1) \cdots {\bm g}(\bm{u}_{N_T}) ]\in\mathbb{R}^{N\times N_T}.
$$
The snapshots matrices satisfy
$$
{\mathcal G }' = {\bm A}_{\mathcal G}   {\mathcal G },
$$
where the unknown matrix ${\bm A}_{\mathcal G}$ is given as the solution of the minimization problem
\begin{equation*}
\min   \left\| {\mathcal G }' -   {\bm A}_{\mathcal G}{\mathcal G }\right\|_{F}^2,
\end{equation*}
where $\|\cdot\|_F$ denotes the Frobenius norm, ${\bm A}_{\mathcal G} = {\mathcal G}'  {\mathcal G}^{\dag}$,
and $\dag$ is the Moore-Penrose pseudoinverse. The Algorithm~\ref{exactDMD_alg} gives the exact DMD procedure \cite{Tu14} to find the DMD basis matrix $\bm{\Psi}^{\text{DMD}}=[\bm{\Psi}^{\text{DMD}}_1 \cdots  \bm{\Psi}^{\text{DMD}}_{\widetilde{m}}]\in\mathbb{R}^{N\times \widetilde{m}}$.
After calculating the DMD basis function $\bm{\Psi}^{\text{DMD}}$ by the Algorithm~\ref{exactDMD_alg},  the nonlinear vector ${\bm g}(\bm{u})$ in the full order state equations \eqref{full_order_ode_state} is approximated as the following time-dependent vector
\begin{equation}\label{dmdapprox}
 \bm{g}^{\text{DMD}}(t) = \sum_{j=1}^{\widetilde{m}}\alpha_j\bm{\Psi}_j^{\text{DMD}} \exp(\omega_ j t) = \bm{\Psi}^{\text{DMD}}\text{diag}(e^{\omega t})\alpha ,
 \end{equation}
 where $\alpha=[\alpha_1,\ldots, \alpha_{\widetilde{m}}]^T$  is the initial amplitudes given by $\alpha = (\bm{\Psi}^{\text{DMD}})^{\dag}{\bm g}(\bm{u}_1)$,  and $\omega =[\omega_1,\ldots, \omega_{\widetilde{m}}]^T$ includes the eigenvalues $\lambda_j$ as $\omega_j=\log{(\lambda_j)}/\Delta t, j=1,\ldots \widetilde{m}$.

Inserting the identity \eqref{dmdapprox} into the reduced state equations \eqref{staterom}, we obtain the following linear ROM for the activator equation
\begin{equation}\label{statedmd}
 \frac{d}{dt} \bm{u}^r +  \bm{S}^r_u \bm{u}^r + \bm{\Psi}_u^T \bm{g}^{\text{DMD}}(t)
  + \bm{M}_{u,v}^r \bm{v}^r  =  \bm{\Psi}_u^T \bm{\ell}_u + \bm{M}_{u,f}^r \bm{f}^r.
\end{equation}
Although the dimension of the reduced system \eqref{statedmd} with $\widetilde{m}$  DMD modes is  small as the reduced system \eqref{statedeim} with $m$ DEIM modes, the POD--DMD reduced state equation \eqref{statedmd} is linear, and the OCP problem becomes convex and Newton's iterations are no more needed.  Therefore, the POD--DMD is much faster than the POD and the POD--DEIM.

\section{Numerical results}
\label{numeric}

We consider the OCP \eqref{ocp}-\eqref{state} with a fast wave speed $V_{max}=64$, and with the parameters and initial condition given by
\[
c_1=9, \quad c_2 = 0.02, \quad c_3 = 5, \quad \epsilon = 0.1, \quad d_{y}=d_{z}=1, \quad \upsilon =10^{-3},
\]
\[
u_0(x,t) = \left\{
             \begin{array}{ll}
               0.1, & \hbox{if} \;\; 0 \leq x_1 \leq 0.1,\; 0 \leq x_2 \leq H  \\
               0, & \hbox{otherwise},
             \end{array}
           \right.
, \quad
v_0(x,0) = 0.
\]
The desired terminal states are chosen as
\[
u_T(x) = u_{\mathrm{nat}}(x,T/2) \quad \hbox{and} \quad v_T(x) = v_{\mathrm{nat}}(x,T/2),
\]
where $u_{\mathrm{nat}}$ and $v_{\mathrm{nat}}$ stand for the solutions of the uncontrolled convective FHN equations. The box constraint for the control is given by
\[
\mathcal{F}_{ad} := \{ f \in L^{\infty}(Q): \; -0.2 \leq f(x,t) \leq 0.2 \;\; \hbox{for   a.e   } (x,t) \in Q \}.
\]
The space domain is taken as a rectangle with $H=5$ and $L=100$. We use uniform step size in space, $\Delta x_{1}=\Delta x_{2} =0.5$, and in time, $\Delta t=0.05$ for the final time $T=1$. Stopping criteria for the OCP is the relative error $|J_{old}-J|/ |J_{old}| \leq 10^{-3}$ on the objective functional.
For $w\in\{u,v,f\}$, the number of POD modes are determined according to the following energy criteria (relative information content (RIC) )
\begin{equation*}
\epsilon(k_w)=\frac{\sum_{i=1}^{k_w}\sigma_{w,i}^2}{\sum_{i=1}^{d_w}\sigma_{w,i}^2},
\end{equation*}
which represents the energy captured by the first $k_w$ POD modes over all $d_w$ POD modes, $d_w$ is the rank of the related snapshot matrix, and $\sigma_{w,i}$ is the corresponding singular value of the $i$-th mode. In the computations, the number of POD modes are taken as the same, $k:=\max\{k_u,k_v,k_f\}$, so that the energy criteria $\min_{k_w} \epsilon(k_w) \geq \%99.99$ is satisfied for $w\in\{u,v,f\}$. Moreover, for $w\in\{u,v,f\}$, we define the following relative $L^2$-error
\[
\|w_h-w_h^r\|:= \frac{\|w_h(x,t)-w_h^r(x,t)\|_{L^2(\Omega)}}{\|w_h(x,t)\|_{L^2(\Omega)}},
\]
to compute the errors between full and reduced order solutions at a given time $t$.

\begin{figure}[H]
\centering
   \includegraphics[scale=0.45]{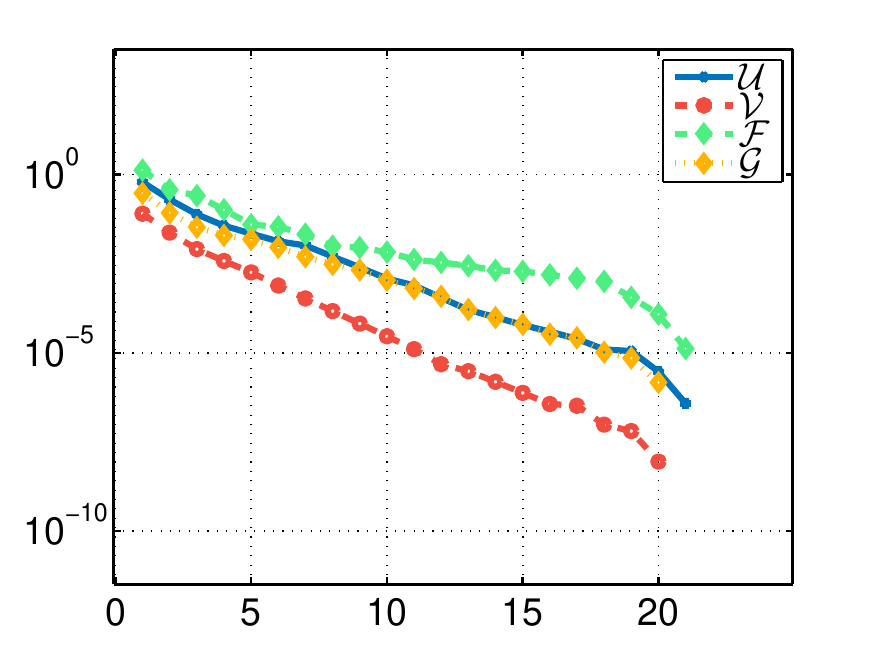}
   \caption{Singular values of the snapshot matrices}
   \label{fig:sing}
\end{figure}

The singular values of the snapshot matrices of the states $u$, $v$, of the control $f$, and of the non-linearity $g(u)$ in Figure~\ref{fig:sing} decay slowly almost at the same rate, which is typical for convection dominated problems. According to the RIC for $\%99.99$ of energy captured,  we have chosen $k=9$ POD modes for the states $u$, $v$ and the control $f$.

\begin{figure}[htb]
\centering
   \subfloat{\includegraphics[scale=0.33]{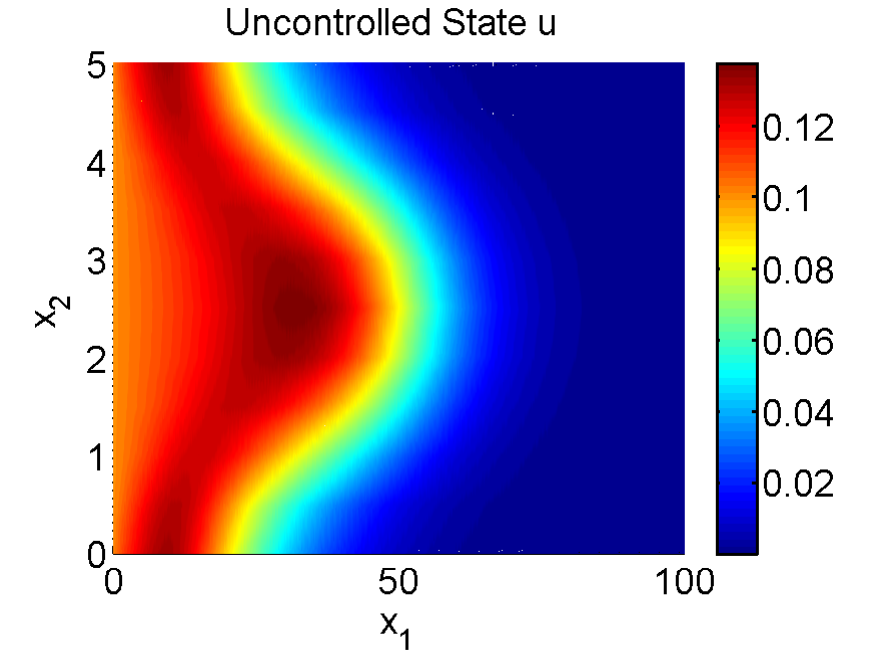}}
   \subfloat{\includegraphics[scale=0.33]{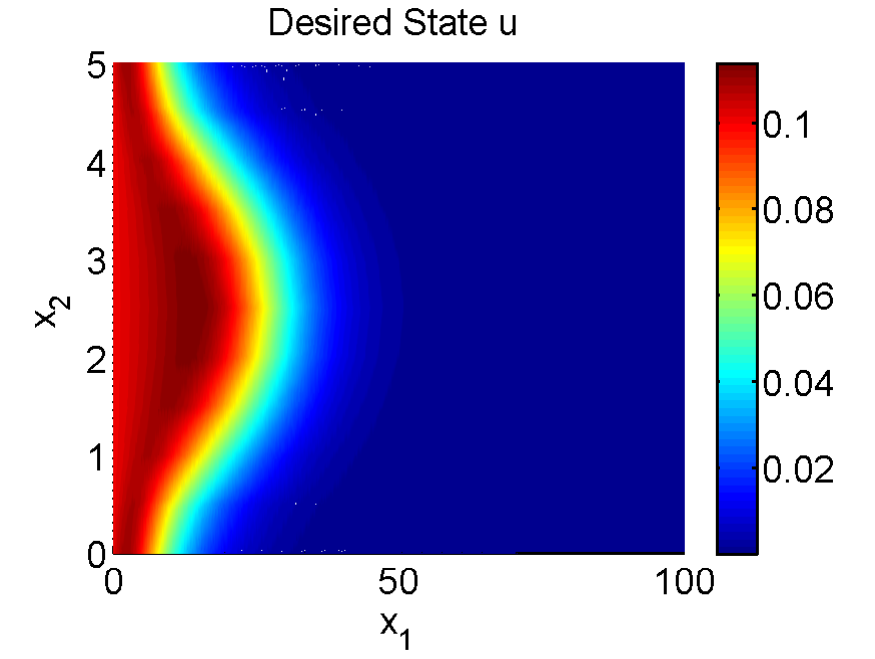}}
   \subfloat{\includegraphics[scale=0.33]{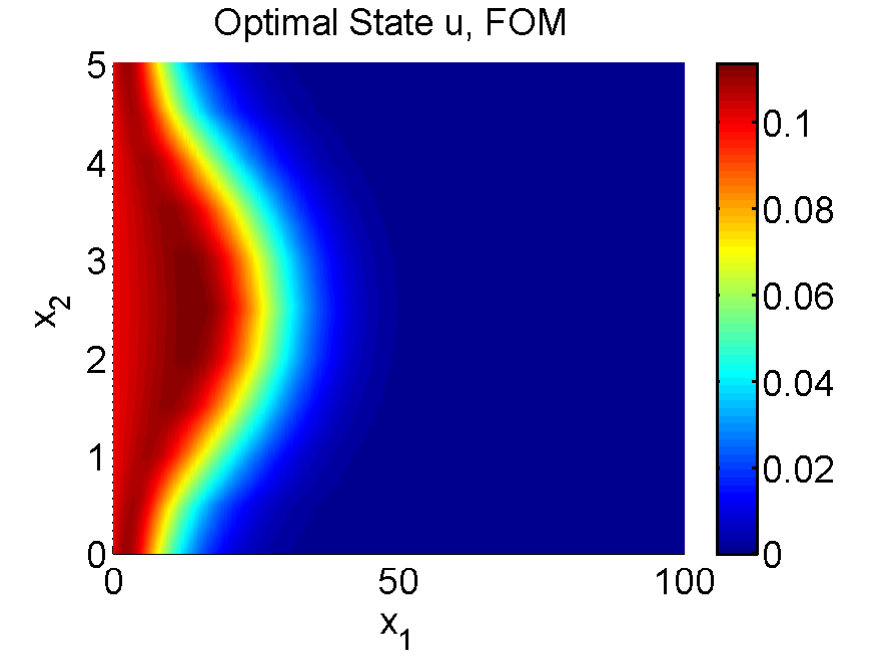}}

   \subfloat{\includegraphics[scale=0.33]{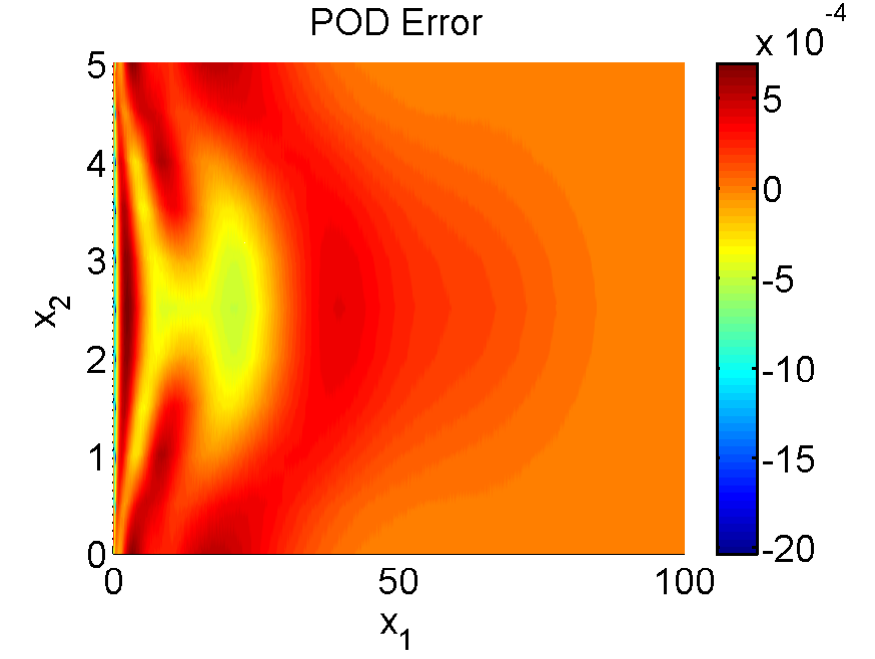}}
   \subfloat{\includegraphics[scale=0.33]{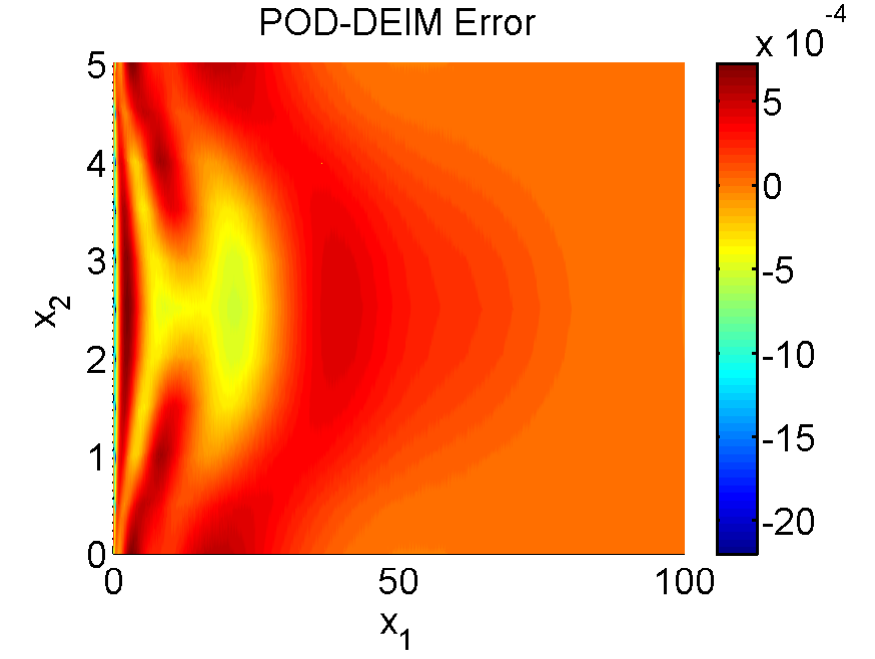}}
	 \subfloat{\includegraphics[scale=0.33]{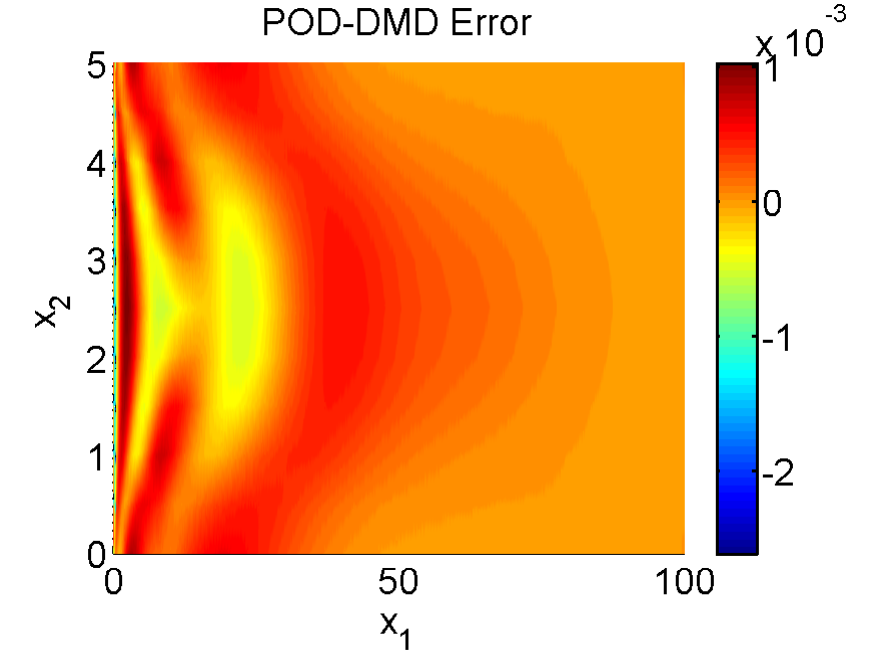}}
   \caption{State u profiles at the final time: (Top) Uncontrolled, desired and optimal FOM states; (Bottom) FOM-ROM errors for $k=9$ POD modes, $m=14$ DEIM basis and $\widetilde{m}=18$ DMD basis.}
   \label{fig:states_u}
\end{figure}

\begin{figure}[htb]
\centering
   \subfloat{\includegraphics[scale=0.33]{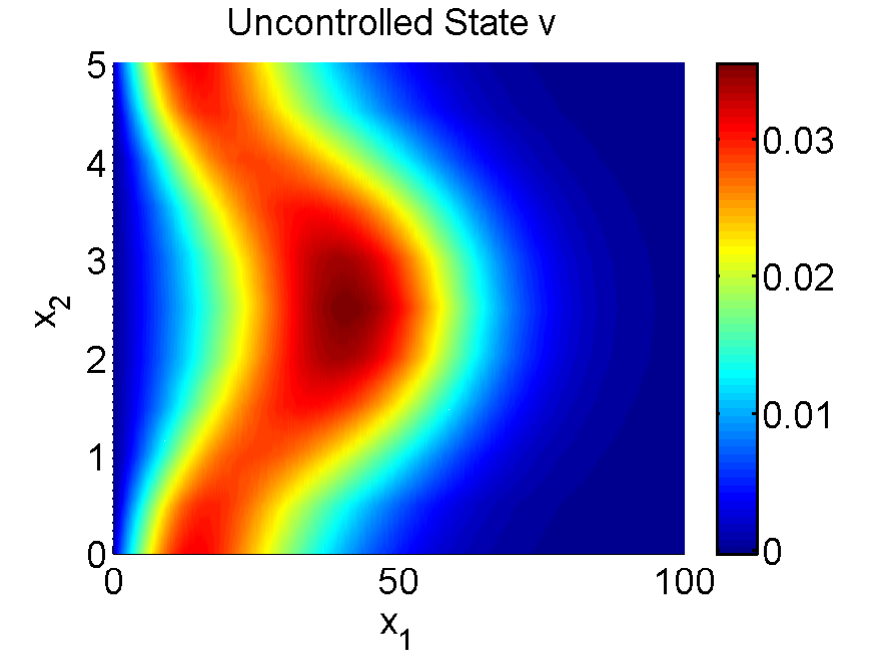}}
   \subfloat{\includegraphics[scale=0.33]{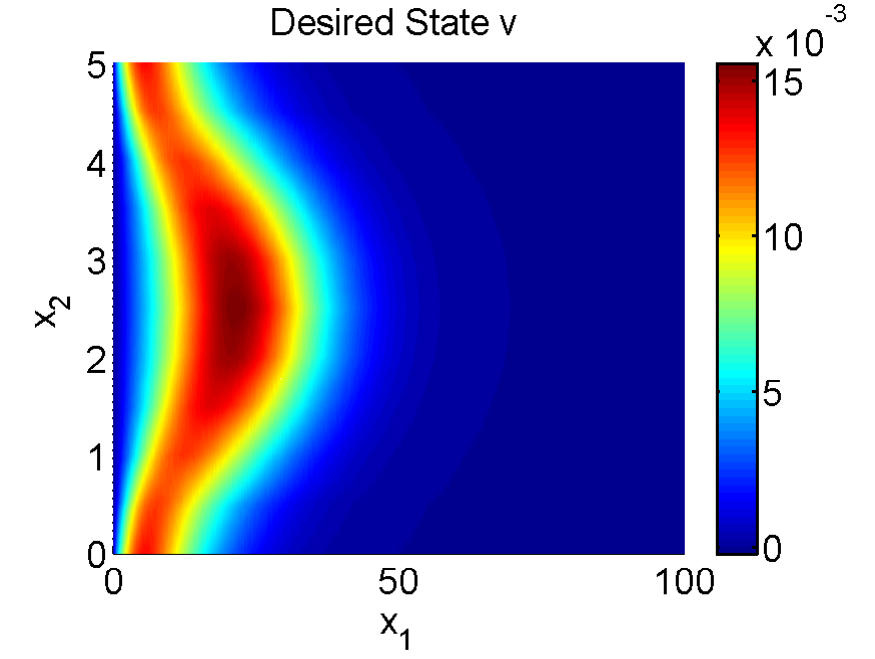}}
   \subfloat{\includegraphics[scale=0.33]{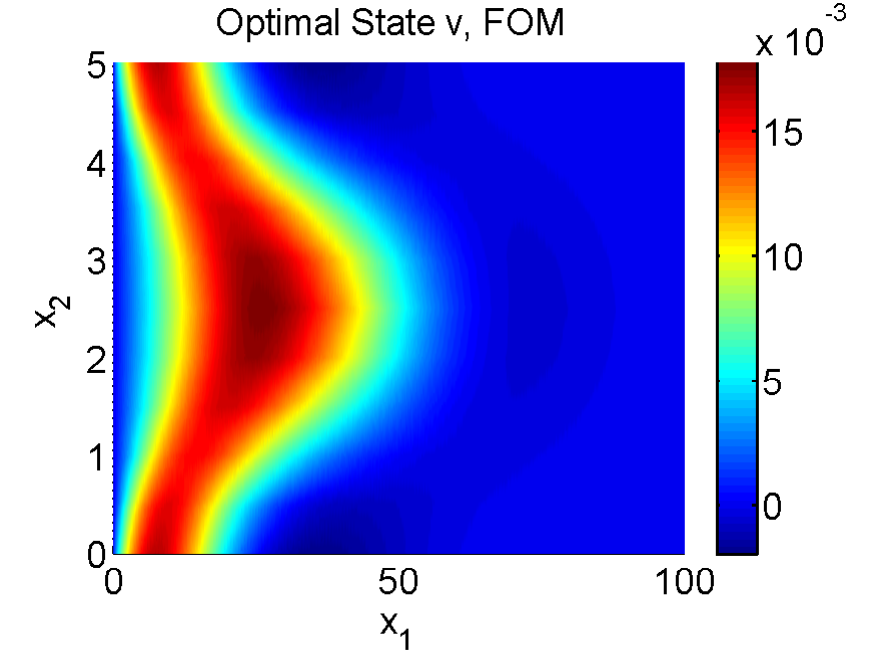}}

   \subfloat{\includegraphics[scale=0.33]{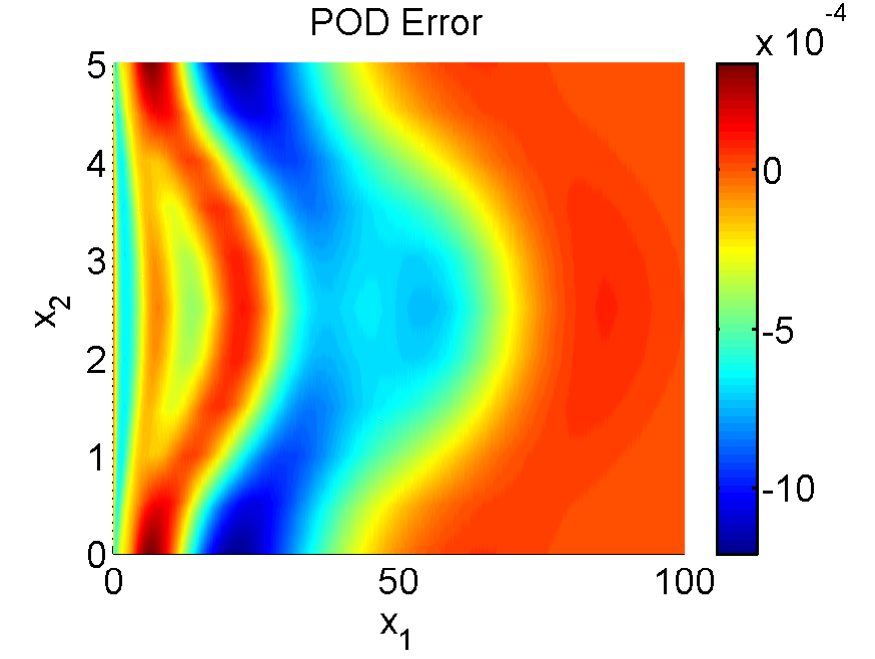}}
   \subfloat{\includegraphics[scale=0.33]{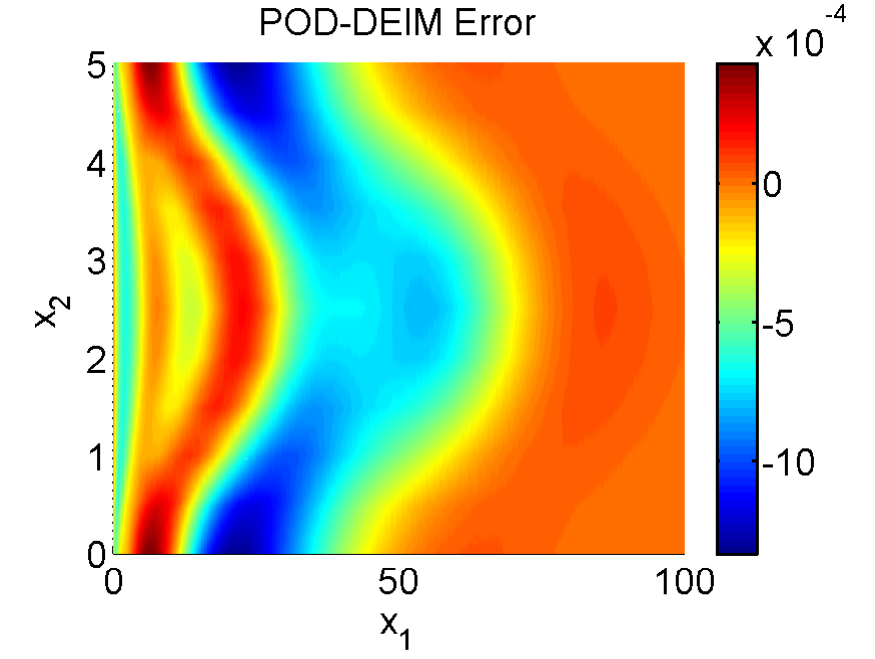}}
	 \subfloat{\includegraphics[scale=0.33]{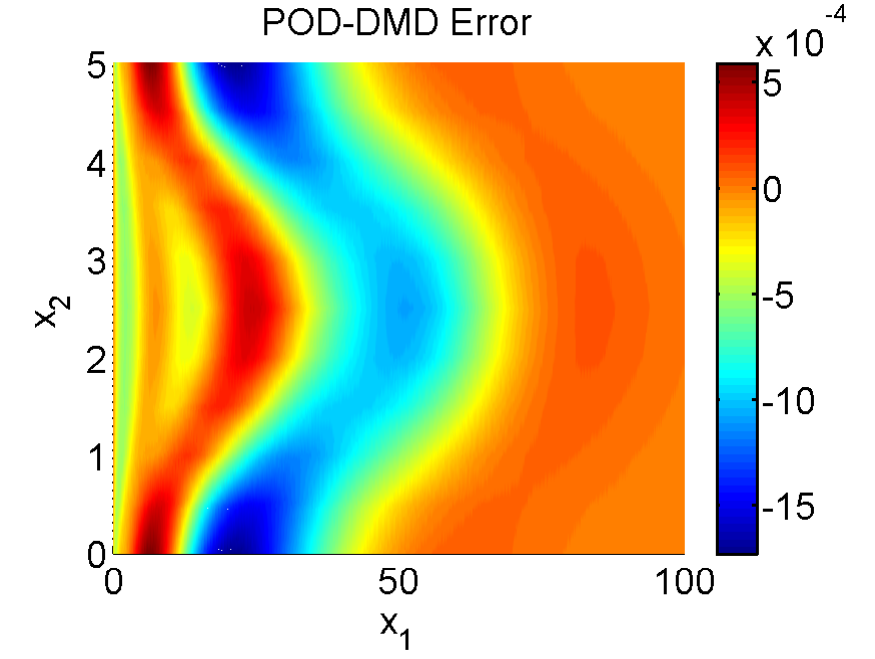}}
   \caption{State v profiles at the final time: (Top) Uncontrolled, desired and optimal FOM states; (Bottom) FOM-ROM errors for $k=9$ POD modes, $m=14$ DEIM basis and $\widetilde{m}=18$ DMD-basis.}
   \label{fig:states_v}
\end{figure}

\begin{figure}[htb]
\centering
   \subfloat{\includegraphics[scale=0.35]{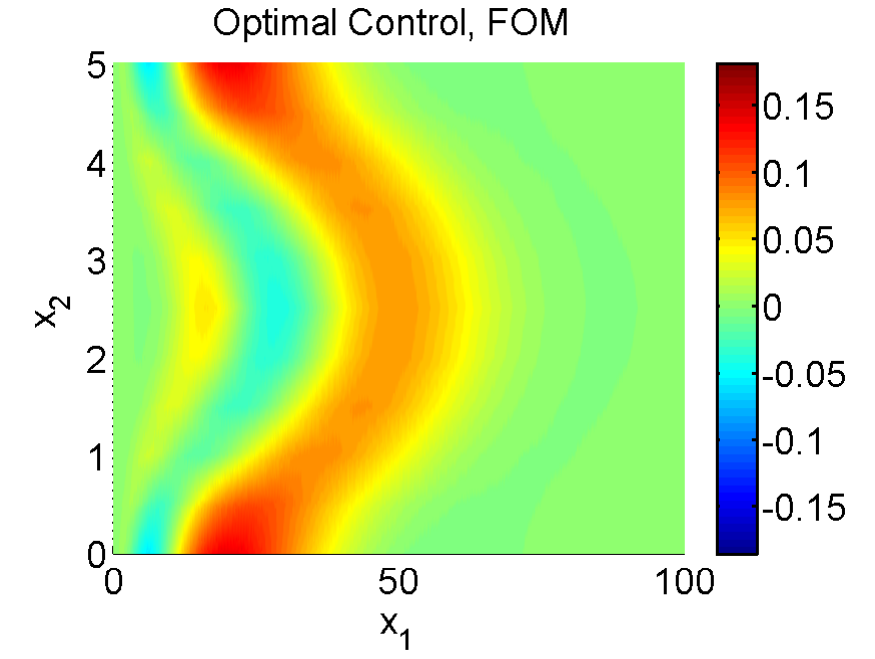}}
   \subfloat{\includegraphics[scale=0.35]{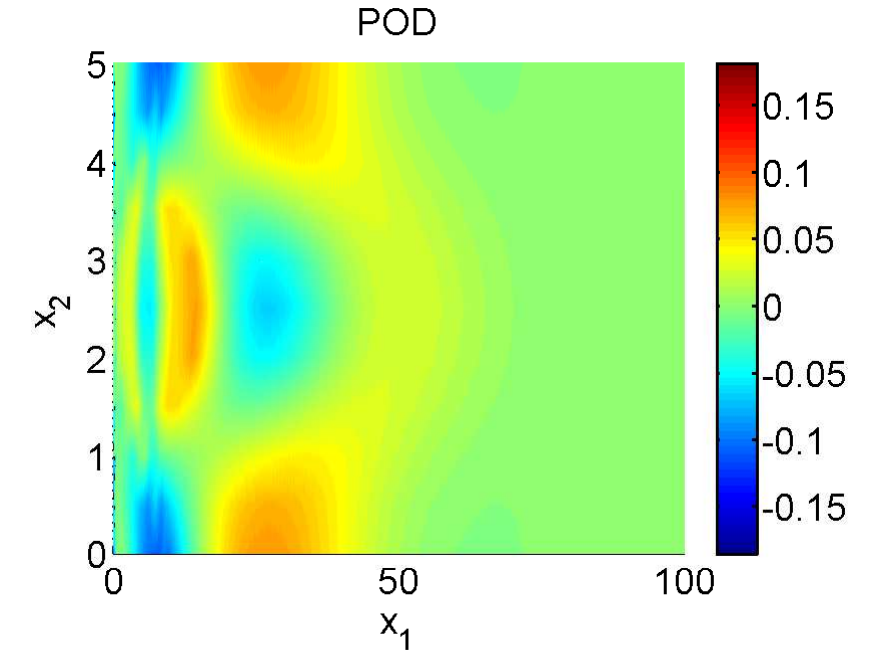}}	
	
	 \subfloat{\includegraphics[scale=0.35]{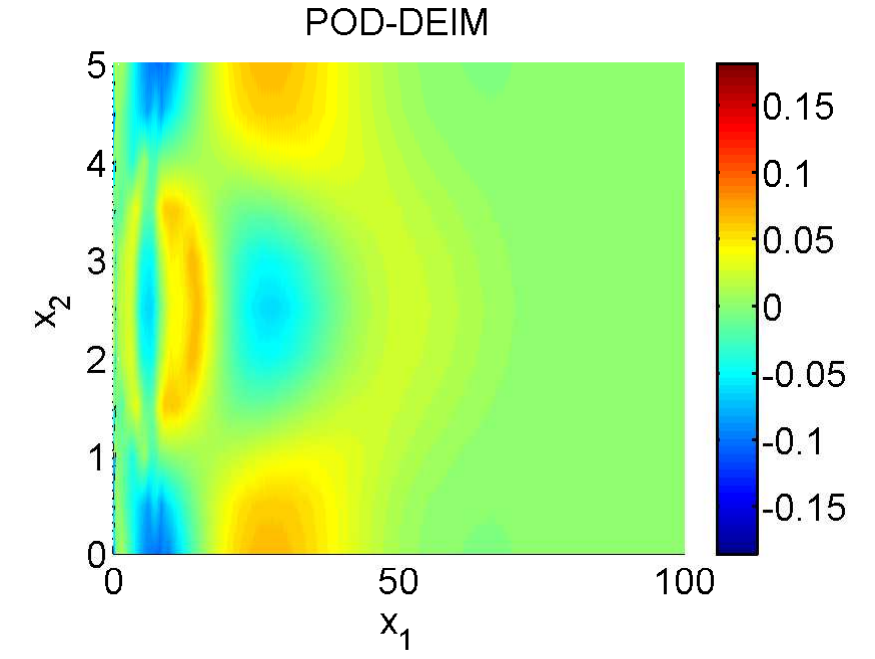}}
   \subfloat{\includegraphics[scale=0.35]{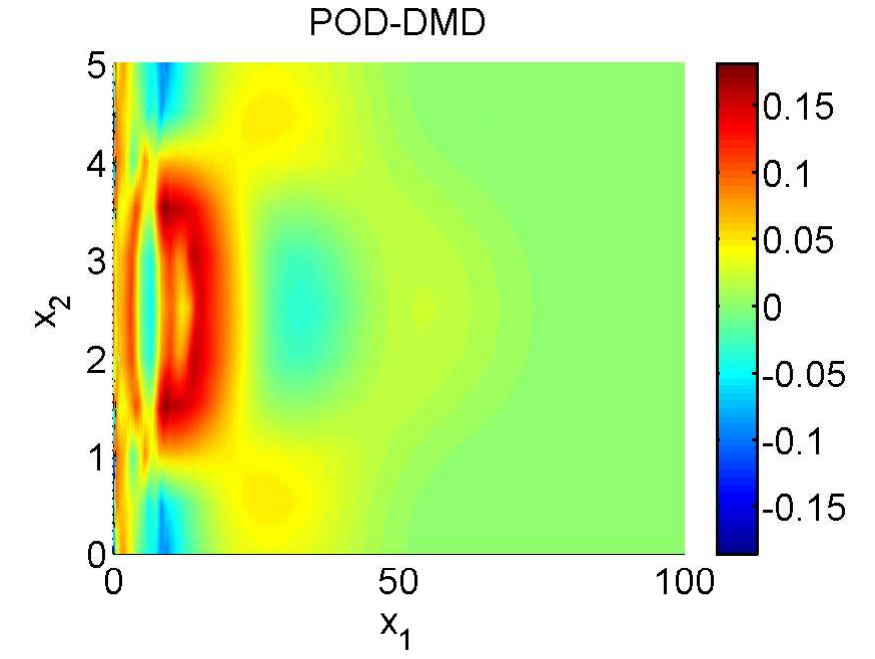}}
   \caption{Control profiles at the final time for $k=9$ POD modes $k=9$, $m=14$ DEIM basis and $\widetilde{m}=18$ DMD basis.}
   \label{fig:control}
\end{figure}

The FOM state solutions at the final time in Figures~\ref{fig:states_u}-\ref{fig:states_v}, top, show the same wave type characteristics  as in \cite{Uzunca17}. The error plots between full and reduced order state solutions in Figures~\ref{fig:states_u}-\ref{fig:states_v}, bottom, agree for POD, POD--DEIM and POD--DMD with $m=14$ and $\widetilde{m}=18$ DEIM and DMD modes, respectively. In Figure~\ref{fig:control}, we give the full and reduced order control profiles at the final time. We see that there is a bit difference between the full and reduced order profiles, but the reduced order profiles are similar in between. The relative $L^2$-errors between the full and reduced order solutions of the states and the control at the final time are shown in Table~\ref{table:time}. Both states $u$ and $v$ are approximated by the reduced models with the same accuracy, whereas the errors in control are larger for POD--DMD due to linearization of the reduced OCP.

\begin{table}[htb]
\centering
\begin{tabular}{l|ccc}
  &  $\|u_h-u_h^r\|$  &  $\|v_h-v_h^r\|$ & $\|f_h-f_h^r\|$  \\
\hline
POD       &  4.644e-03  &  6.167e-02  &  6.614e-01   \\
POD--DEIM  &  4.988e-03  &  6.511e-02  &  7.123e-01  \\
POD--DMD   &  5.787e-03  &  7.670e-02  &  1.061e+00  \\
\hline
\end{tabular}
\caption{Relative $L^2$-errors between full and reduced order states and control at the final time.}
\label{table:time}
\end{table}

\begin{figure}[htb]
\centering
   \subfloat{\includegraphics[scale=0.3]{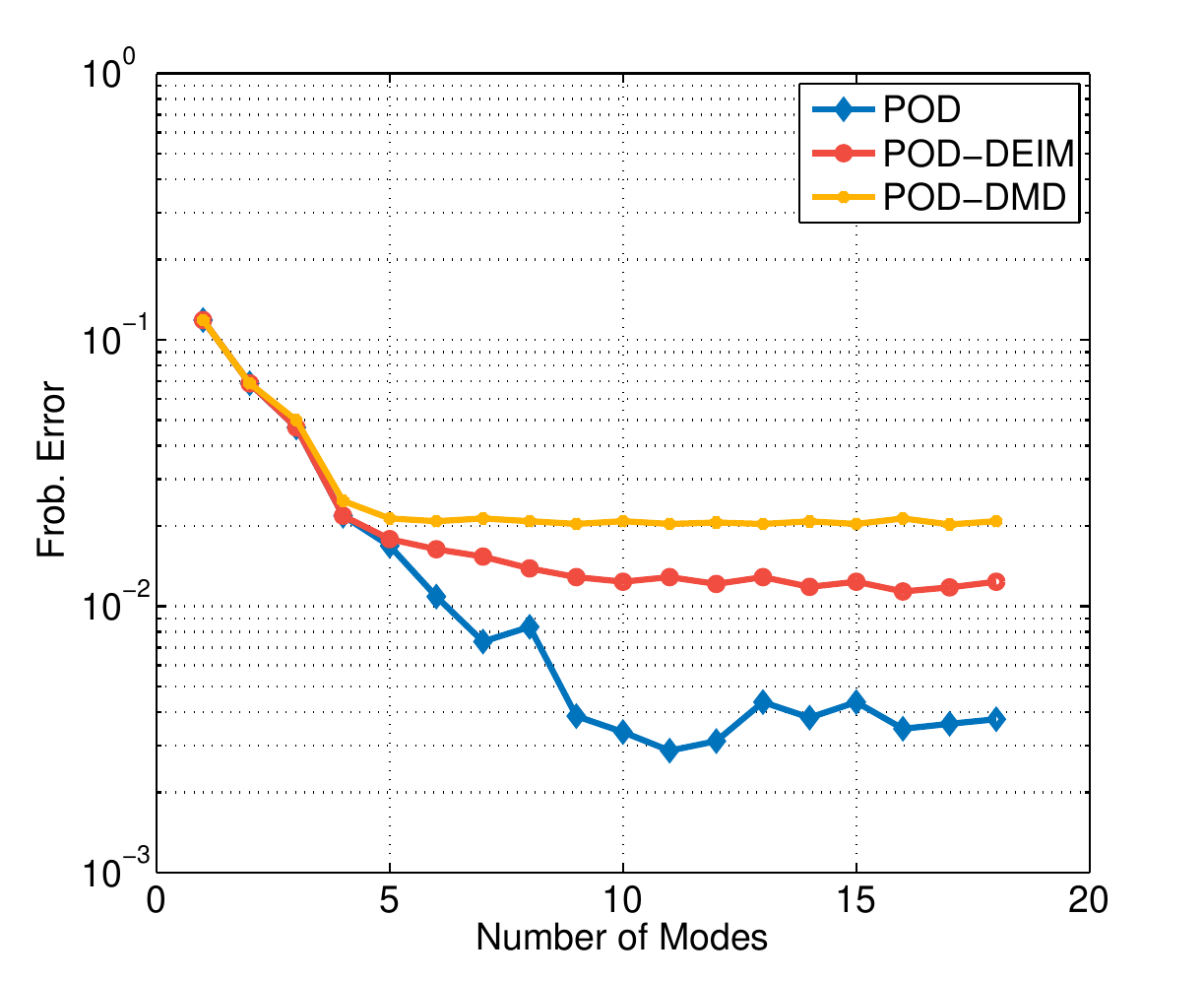}}
   \subfloat{\includegraphics[scale=0.3]{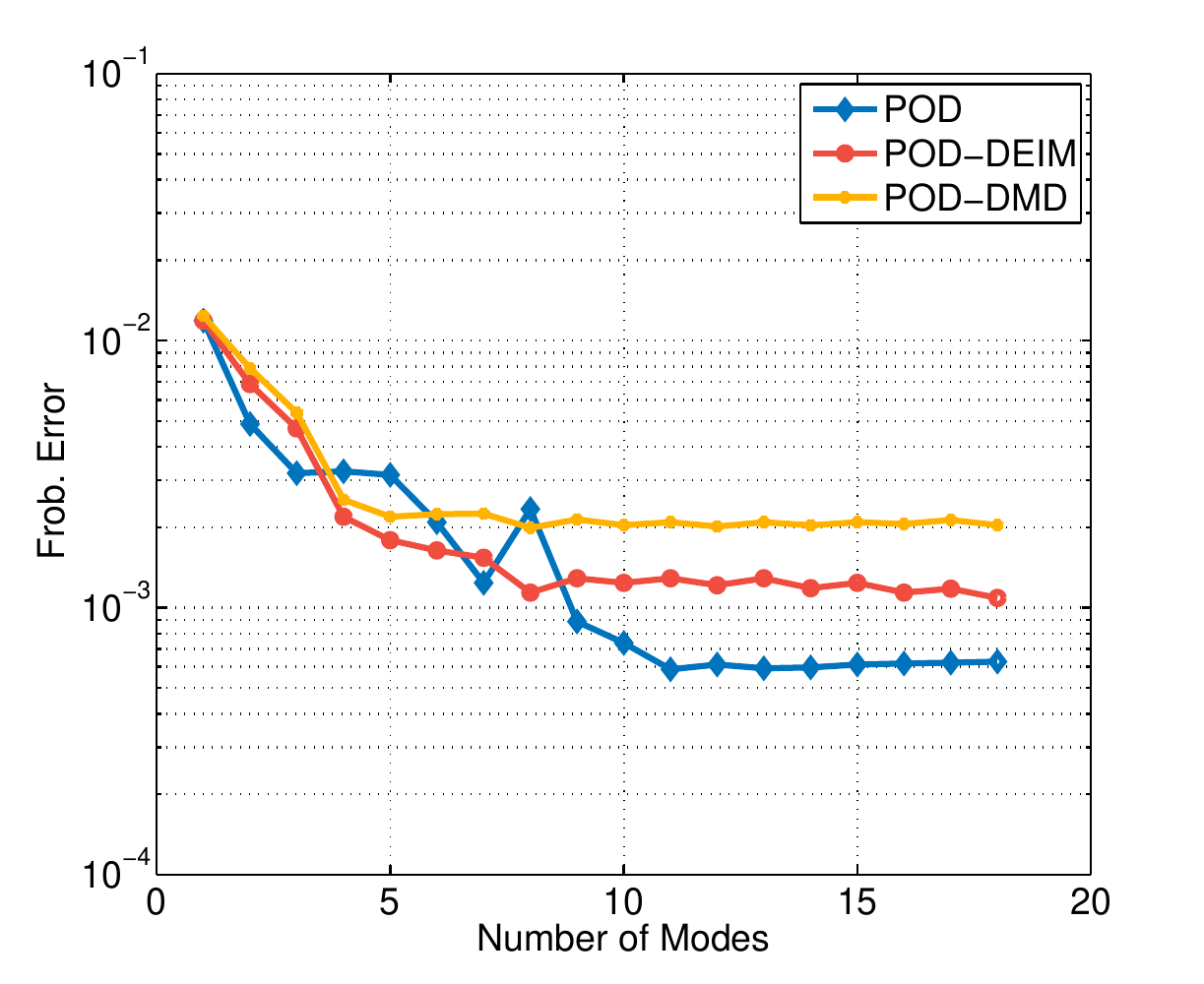}}
   \caption{Relative FOM-ROM Frobenius errors for state $u$ (left) and  state $v$ (right).}
   \label{fig:errors}
\end{figure}

\begin{figure}[htb]
\centering
   \subfloat{\includegraphics[scale=0.3]{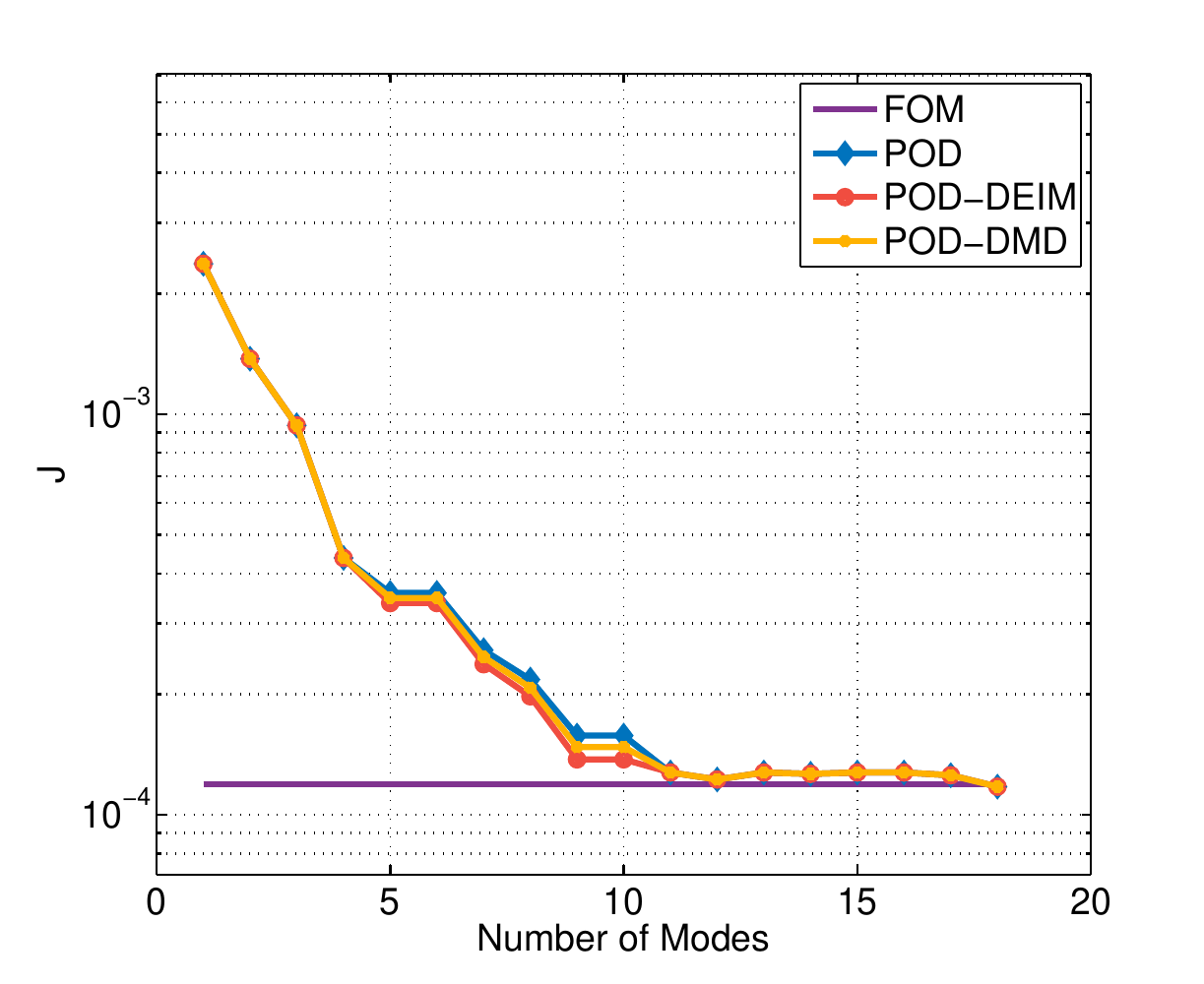}}
   \subfloat{\includegraphics[scale=0.3]{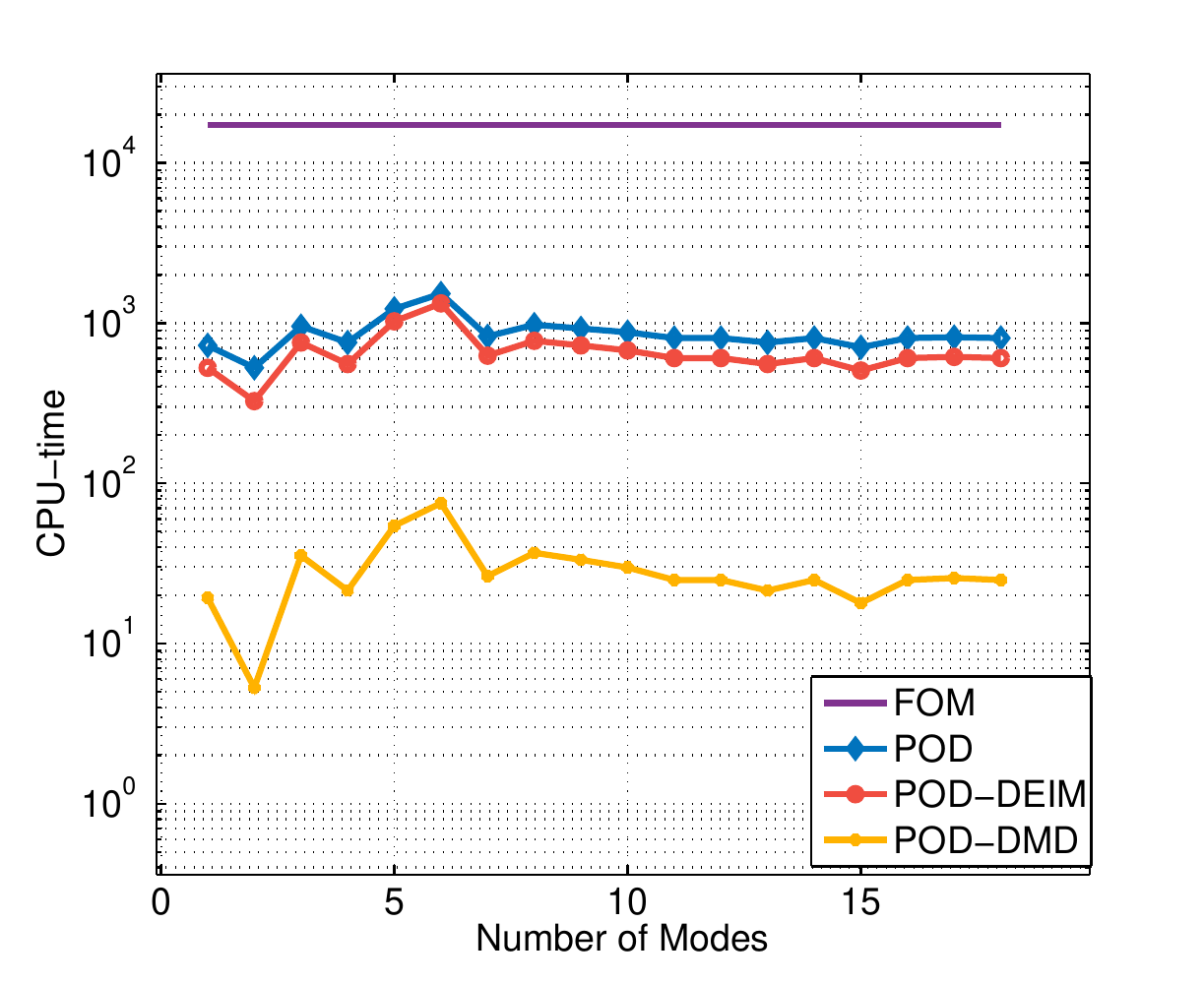}}
   \caption{The values of the objective functionals (left) and CPU times (right) with increasing number of POD modes.}
   \label{fig:cpu}
\end{figure}

Because the Galerkin projected POD--DMD system is linear, the resulting OCP is convex. Therefore, no Newton iteration is recorded for the POD--DMD in Table~\ref{table:all}, and instead of nonlinear CG, we use linear CG method. Comparing the speedup factors (the ratio of the CPU times (in seconds) required for evaluation of FOM to the CPU times required for evaluation of ROM) in Table~\ref{table:all}, the computational efficiency of the POD--DMD approach is clearly visible.   The POD reduced sub-optimal solutions are the most accurate, the POD--DEIM and POD--DMD solutions are less accurate, whereas the POD solutions oscillate at higher modes, as shown for the relative FOM-ROM errors in the Frobenius norm in Figure~\ref{fig:errors}, and in Table~\ref{table:time} as well . The value of the objective function $J$ obtained by full and reduced order solutions are almost the same with increasing number of POD modes in Figure~\ref{fig:cpu}, left.
The ROMs require significantly less CPU time compared to the FOM.
In terms of computational cost, the POD--DMD is the fastest, Figure~\ref{fig:cpu}, right. For a small number of POD modes, the CPU times oscillate until $k \approx 6$, and after they do not change much. Therefore it is not necessary to use a larger number of POD modes considering the CPU times, and the accuracy of the reduced order solutions and reduced order objective functional.

\begin{table}[htb]
\centering
\begin{tabular}{l|c|c|c|c|c|c}
                    & $J$       & $\#$ CG   & $\#$ Line   &  $\#$Newton    & CPU   &  Speedup    \\
                    &         & iterations  &  searches  &   iterations   &  times  &   factors   \\
\hline
  $\text{FOM }$     & 2.376e-03    &  184   &  710   &  3.22   & 15496.7  &    -     \\
  $\text{POD}$      & 2.396e-03    &  85    &  317   &  3.19   & 828.3    &   18.7  \\
  $\text{POD--DEIM}$ & 2.392e-03    &  88    &  332   &  3.20   & 386.2    &   40.1  \\
  $\text{POD--DMD}$  & 2.396e-03    &  99    &  370   &   1     & 8.0      &   1937.1 \\
\hline
\end{tabular}
\caption{The values of the optimal objective functionals $J$, nonlinear CG iterations and line searches, average number of Newton iterations per time step, CPU times and speedup factors.}
\label{table:all}
\end{table}

\clearpage
\section{Conclusions}
\label{conc}

In this study, the computational efficiency and accuracy of three ROM methods are investigated for solving semilinear convection dominated OCP problem. We derived the ROMs for the state and adjoint equations using the same POD basis functions generated by the snapshot matrices of the states only. Two hyperreduction techniques, the DEIM and DMD are applied to reduce the computational cost arising from the nonlinear term in the activator equation. Among the three ROM techniques the POD without DEIM/DMD is the most accurate as expected. The POD--DEIM and POD--DMD errors are close, but the POD--DMD is the fastest due to the fact that the reduced state equations are no more nonlinear and the OCP is convex. In a future study, we plan to compare these ROM techniques for model predictive control and feedback control problems using compressive POD and DMD.\\

\noindent {\bf Acknowledgements}
The authors would like to thank the reviewers for the comments and suggestions that
helped to  improve the manuscript.


\end{document}